\documentclass[review,12pt,3p]{elsarticle}

\usepackage{tikz}
\usepackage{amsmath} 
\usepackage{amssymb}
\usetikzlibrary{shapes,arrows,chains,intersections}
\usepackage{subfigure}
\usepackage{float}
\usepackage{textcomp} 
\usepackage{url}

\graphicspath{{./pdf/}}

\usepackage{natbib}
\usepackage{multicol} 
\usepackage{nomencl}
\usepackage{framed} 
\makenomenclature
\renewcommand*\nompreamble{\begin{multicols}{2}}
\renewcommand*\nompostamble{\end{multicols}}

\newtheorem{remark}{Remark}

\usepackage{multirow}

\newcommand{\vect}[1]{\boldsymbol{#1}}

\begin{document}

\begin{frontmatter}

\title{Reformulated dissipation for the free-stream preserving of the conservative finite difference schemes on curvilinear grids}

\author[1]{Hongmin Su}
\ead{hongminsu@mail.nwpu.edu.cn}

\author[1]{Jinsheng Cai}
\ead{caijsh@nwpu.edu.cn}

\author[1]{Shucheng Pan\corref{cor1}}  
\ead{shucheng.pan@nwpu.edu.cn}

\author[2]{Xiangyu Hu}
\ead{xiangyu.hu@tum.de}

\address[1]{Department of Fluid Mechanics, School of Aeronautics, Northwestern Polytechnical University, 710072 Xi'an , China}
\address[2]{Department of Mechanical Engineering, Technical University of Munich, 85748 Garching, Germany}

\cortext[cor1]{Corresponding author}

\begin{abstract}
In this paper, we develop a new free-stream preserving (FP) method for high-order upwind conservative finite-difference (FD) schemes on the curvilinear grids. 
This FP method is constrcuted by subtracting a reference cell-face flow state from each cell-center value in the local stencil of the original upwind conservative FD schemes, which effectively leads to a reformulated dissipation. It is convenient to implement this method, as it does not require to modify the original forms of the upwind schemes. In addition, the proposed method removes the constraint in the traditional FP conservative FD schemes that require a consistent discretization of the mesh metrics and the fluxes. With this, the proposed method is more flexible in simulating the engineering problems which usually require a low-order scheme for their low-quality mesh, while the high-order schemes can be applied to approximate the flow states to improve the resolution. After demonstrating the strict FP property and the order of accuracy by two simple test cases, we consider various validation cases, including the supersonic flow around the cylinder, the subsonic flow past the three-element airfoil, and the transonic flow around the ONERA M6 wing, etc., to show that the method is suitable for a wide range of fluid dynamic problems containing complex geometries.
Moreover, these test cases also indicate that the discretization order of the metrics have no significant influences on the numerical results if the mesh resolution is not sufficiently large.
\end{abstract}


\end{frontmatter}

\section{Introduction}
Recently, high-order numerical methods have become a power tool to understand the underlying mechanism in many complicated fluid mechanics problems, e.g. hydrodynamic instabilities, turbulence, and aeroacoustics. 
Among those methods, the finite difference method (FDM) is the most efficient one, especially in multiple dimensions, as the discretization is performed in a dimension-by-dimension manner, rather than the multi-dimensional operations of e.g. the finite volume method (FVM) and finite element method (FEM) ~\cite{shu2003high}. 
Thus, numerous high-order FD schemes have been developed, such as the compact difference scheme \cite{lele1992compact,VISBAL2002155}, the essentially non-oscillatory scheme (ENO) \cite{HARTEN19973}, the weighted essentially non-oscillatory scheme (WENO) \cite{jiang_efficient_1996}, and the weighted compact nonlinear scheme (WCNS) \cite{DENG200022}. 
However, when the high-order conservative FD schemes are applied on the body-fitted mesh for flows around complex geometries of e.g. aircraft and automobile, the violation of the geometric conservation law (GCL)~\cite{thomas1979geometric}, including the volume conservation law (VCL) ~\cite{ABE201314, abe2014geometric} and the surface conservation law (SCL) ~\cite{zhang1993discrete}, leads to large errors that may destabilize numerical simulations or induces spurious hydrodynamic fluctuations. Specifically for stationary curvilinear grid, this is reflected by the free-stream preservation problem, i.e. the initial uniform flow field becomes nonuniform.

Different numerical treatments have been proposed to preserve the free-stream condition for FDM in general curvilinear grids, such as the pioneer work of Thomas and Lombard \cite{thomas1979geometric} who used the conservative form of metrics for computing geometric transformation. Latter, this conservative metric technique was applied by Visbal and Gaitonde~\cite{VISBAL2002155} to maintain GCL for high-order compact central schemes. A milestone in this area, called the symmetrical conservative metric method (SCMM), was proposed by Deng et al.~\cite{DENG201390} and Abe et al.~\cite{abe2014geometric} independtly, based on the symmetric forms of metrics derived by Vinokur and Yee \cite{vinokur2002extension}. In such way, the asymmetric metric errors of the previously developed conservative metric methods \cite{thomas1979geometric, DENG20111100} is eliminated by preserving the coordinate-invariant property. Other variants of the SCMM are also derived in the literature, such as Liao's cell-center version SCMM \cite{Fei2015Extending1,Fei2015Extending2}. And Abe et. al~\cite{abe2014geometric} elaborate the geometric interpretations of metrics and Jacobian for both the CMM and SCMM. Once this SCMM is applied, a sufficient condition to preserve GCL is to apply the identical central finite difference discretization schemes for the derivatives in both the fluxes and the grid metrics, as proved by Deng and Abe ~\cite{ DENG201390,abe2014geometric}.

Straightforward extension of the above mentioned to upwind linear schemes or nonlinear schemes can not preseve the GCL as the different numerical schemes are requied by grid metrics and fluxes, which violate the sufficient condition given by Deng and Abe ~\cite{DENG201390,abe2014geometric}. One way to address this is to independently approximate the flow variables and grid metrics. A successful strategy to maintain GCL for nonlinear schemes is the combination of the conservative metrics and the WCNS, as demonstrated by Nonomura et al.~\cite{Nonomura2010_CF}. Similar method is the alternative finite-difference form of WENO (AWENO)~\cite{JiangAn, jiang2018, Yu2020}. Another way is to deal with the central part and the dissipative part of the upwind scheme, respectively. For example, employing the high-order central schemes to the former and performing the finite-volume-like schemes, i.e. freezing metrics either for the entire stencil ~\cite{nonomura_new_2015}, or for the local difference form partially \cite{zhu_free-stream_2019}, or replacing the transformed conservative variables with the original one ~\cite{LiFurther}, to the dissipative part. In addition, to satisfy the above sufficient condition, Ref. \cite{2017Zhu} introduces an offsetting terms with the same WENO nonlinear weights for computing the corresponding inviscid fluxes. Nevertheless, further improvement is required when upwind linear schemes or nonlinear schemes are applied for realistic engineering problems. On the one hand, these methods modify the standard forms of the linear upwind and WENO schemes, which yields operation complexity and difficulty for the implementation. On the other hand, these methods are required to satisfy the sufficient condition provided by Deng and Abe ~\cite{DENG201390,abe2014geometric}, which introduces a constraint that they must apply the unique high-order discretization for the grid metrics and the central part of the fluxes. As a result, the simulations may blow up on a low quality grid.

In this study, we propose an efficient strategy to maintain the free-stream preserving identity by applying a reformulated upwind dissipation to the linear upwind and WENO schemes. The novelty is twofold. First, the standard forms of the original schemes are not required to be modified, indicating an easy and starightfoward implementation. Second, they remove the constraint of using the consistent discretization for the grid metrics and fluxes, which permits using low-order schemes for grid information to enhance the robustness in low-quality meshes and using high-order schemes for flow fluxes to increase the effective resolution. The outline is organized as follows. In Sec. \ref{sec:govern}, we introduce the Navier-Stokes (NS) equations, metrics and SCL identity. The introduction of the linear upwind and standard WENO schemes on curvilinear grids based on flux vector splitting is also included in this section. In Sec. \ref{sec:method}, the proposed strategy on free-stream preserving for linear upwind and WENO schemes is described in detail. 
Free-streaming preserving capability, order of accuracy, and robustness of the proposed method are demonstrated in Sec. \ref{sec:testcases} by a range of numerical examples, followed by a concluding remark in Sec. \ref{sec:conlusion}.

\section{Governing equations and numerical methods} \label{sec:govern}
\subsection{Governing equations and metrics on stationary curvilinear coordinates}
The compressible Navier-Stokes equations on curvilinear grids are given by
\begin{equation} \label{eq:NS} 
   \begin{aligned}
\dfrac{\partial \vect{Q}}{\partial t}&+\dfrac{\partial}{\partial \xi} \left( \dfrac{\xi_x \vect{F}+\xi_y \vect{G}+\xi_z \vect{H}}{J} \right)+\dfrac{\partial}{\partial \eta} \left( \dfrac{\eta_x \vect{F}+\eta_y \vect{G}+\eta_z \vect{H}}{J} \right)+\dfrac{\partial}{\partial \zeta} \left( \dfrac{\zeta_x \vect{F}+\zeta_y \vect{G}+\zeta_z \vect{H}}{J} \right)\\
& -\dfrac{\partial}{\partial \xi} \left( \dfrac{\xi_x \vect{F_v}+\xi_y \vect{G_v}+\xi_z \vect{H_v}}{J} \right)-\dfrac{\partial}{\partial \eta} \left( \dfrac{\eta_x \vect{F_v}+\eta_y \vect{G_v}+\eta_z \vect{H_v}}{J} \right)-\dfrac{\partial}{\partial \zeta} \left( \dfrac{\zeta_x \vect{F_v}+\zeta_y \vect{G_v}+\zeta_z \vect{H_v}}{J} \right)=\vect{0}
\end{aligned}
\end{equation}
with
\begin{equation} \label{eq:Q} 
  \boldsymbol{Q}=\left(
        \begin{matrix}
          \rho & \rho u_1  &  \rho u_2   &  \rho u_3    &   \rho E 
        \end{matrix}   
        \right)^T,
\end{equation}
\begin{equation} \label{eq:F}
\boldsymbol{F}=\left(
\begin{matrix}
          \rho u_1 & \rho u_1 u_1+p & \rho u_2 u_1 & \rho u_3 u_1  & (\rho E+p) u_1
\end{matrix}
        \right)^T,
\end{equation}        
\begin{equation} \label{eq:G}        
\boldsymbol{G}=\left(
\begin{matrix}
          \rho u_2 & \rho u_1 u_2 & \rho u_2 u_2+p  &\rho u_3 u_2 & (\rho E+p) u_2
\end{matrix}
        \right)^T,
\end{equation}  
\begin{equation} \label{eq:H}        
\boldsymbol{H}=\left(
\begin{matrix}
          \rho u_3 & \rho u_1 u_3 & \rho u_2 u_3 & \rho u_3 u_3+p  & (\rho E+p) u_3
\end{matrix}
        \right)^T,
\end{equation} 
\begin{equation} \label{eq:NSv1}
\vect{F_v}=\left(
        \begin{matrix}
           0& \tau_{11} & \tau_{12} & \tau_{13}& u_i\tau_{i1}-\dot{q}_1 
        \end{matrix}
        \right)^T
\end{equation}        
\begin{equation} \label{eq:NSv2}        
\vect{G_v}=\left(
        \begin{matrix}
           0& \tau_{21} & \tau_{22} & \tau_{23}& u_i\tau_{i2}-\dot{q}_2 
        \end{matrix}
        \right)^T
\end{equation}        
\begin{equation} \label{eq:NSv3}        
\vect{H_v}=\left(
        \begin{matrix}
           0& \tau_{31} & \tau_{32} & \tau_{33}& u_i\tau_{i3}-\dot{q}_3 
        \end{matrix}
        \right)^T                
\end{equation}
where $u_1$, $u_2$, $u_3$, $\boldsymbol{F}$, $\boldsymbol{G}$, $\boldsymbol{H}$ and $\boldsymbol{F_v}$, $\boldsymbol{G_v}$, $\boldsymbol{H_v}$ denote the velocity components, the inviscid and viscous flux vectors in x, y and z direction, respectively. $\rho$, $p$ and $E$ are the density, pressure and the total specific energy. t stands for the physical time. $\xi,\eta,\zeta$ are the transformed coordinates on a uniform computational domain, and $J$ is the transformed Jacobian. $\tau_{ij}$ is the shear stress tensor 
\begin{equation} \label{eq:stress}
   \tau_{ij}=2\mu (S_{ij}-\delta_{ij}\dfrac{S_{kk}}{3}),
\end{equation}
\begin{equation} \label{eq:strain} 
   S_{ij}=\dfrac{1}{2}(\dfrac{\partial u_i}{\partial x_j}+\dfrac{\partial u_j}{\partial x_i}),
\end{equation}
and $\dot{q}_i$ is the heat flux in direction i
\begin{equation} 
  \dot{q}_{i}=-\lambda \dfrac{\partial T}{\partial x_i},
\end{equation}
where $\mu$ and $\lambda$ is the shear viscosity and thermal conductivity. 

The equation of state for ideal gas is
\begin{equation} \label{eq:EOS} 
   p=\left( \gamma-1\right) \rho e,
\end{equation}    
where the specific heat ratio is $\gamma=1.4$. For convenience, we denote the fluxes in computational space as
\begin{equation} \label{eq:Q} 
   \boldsymbol{\tilde{Q}}=\dfrac{ \boldsymbol{Q} }{J}
\end{equation}
\begin{equation} \label{eq:FGH} 
\begin{aligned}
\boldsymbol{\tilde{F}}&=\dfrac{\xi_x \boldsymbol{F}+\xi_y \boldsymbol{G}+\xi_z   \boldsymbol{H}}{J} \\
   \boldsymbol{\tilde{G}}&=\dfrac{\eta_x \boldsymbol{F}+\eta_y \boldsymbol{G}+\eta_z \boldsymbol{H}}{J}\\
   \boldsymbol{\tilde{H}}&=\dfrac{\zeta_x \boldsymbol{F}+\zeta_y \boldsymbol{G}+\zeta_z \boldsymbol{H}}{J}
\end{aligned}
\end{equation}
\begin{equation} \label{eq:FGH_v} 
\begin{aligned}
   \boldsymbol{\tilde{F_v}}&=\dfrac{\xi_x \boldsymbol{F_v}+\xi_y \boldsymbol{G_v}+\xi_z \boldsymbol{H_v}}{J}\\
   \boldsymbol{\tilde{G_v}}&=\dfrac{\eta_x \boldsymbol{F_v}+\eta_y \boldsymbol{G_v}+\eta_z \boldsymbol{H_v}}{J}\\
\boldsymbol{\tilde{H_v}}&=\dfrac{\zeta_x \boldsymbol{F_v}+\zeta_y \boldsymbol{G_v}+\zeta_z \boldsymbol{H_v}}{J}    
\end{aligned}
\end{equation}

When subtracting the flow information from the governing equations by applying the free-stream condition to Eq. \eqref{eq:NS}, these equations are simplified as
\begin{equation} \label{eq:metrics_ns} 
\dfrac{\partial}{\partial \xi} \left( \dfrac{\xi_x +\xi_y +\xi_z }{J} \right)+\dfrac{\partial}{\partial \eta} \left( \dfrac{\eta_x +\eta_y +\eta_z }{J} \right)+\dfrac{\partial}{\partial \zeta} \left( \dfrac{\zeta_x +\zeta_y +\zeta_z }{J} \right)
=\vect{0}.
\end{equation}
An sufficient condition to satisfy Eq. \eqref{eq:metrics_ns} is 
\begin{equation} \label{eq:SCL} 
\begin{aligned}
I_x= \left( \dfrac{ \xi_x}{J} \right)_{\xi}
    +\left( \dfrac{ \eta_x}{J} \right)_{\eta}
    +\left( \dfrac{ \zeta_x}{J}\right)_{\zeta}=0, \\
I_y= \left( \dfrac{ \xi_y}{J} \right)_{\xi}
    +\left( \dfrac{ \eta_y}{J} \right)_{\eta}
    +\left( \dfrac{ \zeta_y}{J}\right)_{\zeta}=0,\\
I_z= \left( \dfrac{ \xi_z}{J} \right)_{\xi}
    +\left( \dfrac{ \eta_z}{J} \right)_{\eta}
    +\left( \dfrac{ \zeta_z}{J}\right)_{\zeta}=0.
\end{aligned}
\end{equation}  
These equations represent the consistence of vectorized computational cell surfaces in FVM as discussed in Ref.~\cite{VINOKUR19891}, which are regarded as the SCL proposed by Zhang et al.~\cite{zhang1993discrete}. Although Eq.~\eqref{eq:SCL} are strictly satisfied theoretically, its the numerical discretization errors may violate these identities. To maintian the SCL in high-order central finite difference schemes, the SCMM can be used to discretize the grid metrics by adopting the identical scheme for the fluxes discretization~\cite{DENG201390,abe2014geometric}. The symmetry conservative metrics can be expressed by
\begin{equation} \label{eq:JBabc} 
\begin{aligned}
  \dfrac{\xi_x}{J}=\dfrac{1}{2} \left[   
        \left( y_\eta z \right)_\zeta  
     -  \left( y_\zeta z \right)_\eta 
     +  \left( yz_\zeta\right)_\eta 
     -  \left( yz_\eta \right)_\zeta  
     \right],  \\
  \dfrac{\xi_y}{J}=\dfrac{1}{2} \left[   
        \left( xz_\eta \right)_\zeta  
     -  \left( xz_\zeta \right)_\eta 
     +  \left( x_\zeta z\right)_\eta 
     -  \left( x_\eta z \right)_\zeta  
     \right],  \\  
  \dfrac{\xi_z}{J}=\dfrac{1}{2} \left[   
        \left( x_\eta y \right)_\zeta  
     -  \left( x_\zeta y \right)_\eta 
     +  \left( xy_\zeta\right)_\eta 
     -  \left( xy_\eta \right)_\zeta  
     \right],  \\    
  \dfrac{\eta_x}{J}=\dfrac{1}{2} \left[   
        \left( y_\zeta z \right)_\xi  
     -  \left( y_\xi z \right)_\zeta 
     +  \left( yz_\xi  \right)_\zeta 
     -  \left( yz_\zeta \right)_\xi  
     \right],  \\  
  \dfrac{\eta_y}{J}=\dfrac{1}{2} \left[   
        \left( xz_\zeta\right)_\xi  
     -  \left( xz_\xi  \right)_\zeta 
     +  \left( x_\xi z  \right)_\zeta 
     -  \left( x_\zeta z \right)_\xi  
     \right],  \\      
  \dfrac{\eta_z}{J}=\dfrac{1}{2} \left[   
        \left( x_\zeta y\right)_\xi  
     -  \left( x_\xi y \right)_\zeta 
     +  \left( xy_\xi  \right)_\zeta 
     -  \left( xy_\zeta \right)_\xi  
     \right],  \\  
  \dfrac{\zeta_x}{J}=\dfrac{1}{2} \left[   
        \left(y_\xi z\right)_\eta
     -  \left(y_\eta z \right)_\xi 
     +  \left( yz_\eta  \right)_\xi 
     -  \left( yz_\xi \right)_\eta  
     \right],  \\   
  \dfrac{\zeta_y}{J}=\dfrac{1}{2} \left[   
        \left(xz_\xi \right)_\eta
     -  \left(xz_\eta \right)_\xi 
     +  \left( x_\eta z  \right)_\xi 
     -  \left( x_\xi z \right)_\eta  
     \right],  \\   
  \dfrac{\zeta_z}{J}=\dfrac{1}{2} \left[   
        \left(x_\xi y \right)_\eta
     -  \left(x_\eta y \right)_\xi 
     +  \left( xy_\eta  \right)_\xi 
     -  \left( xy_\xi  \right)_\eta  
     \right],  \\                          
\end{aligned}  
\end{equation}
and 
\begin{equation} \label{eq:JB} 
  \dfrac{1}{J}=\dfrac{1}{3} \left[   
        \left( x \dfrac{\xi_x}{J}+y \dfrac{\xi_y}{J}+z \dfrac{\xi_z}{J}  \right)_\xi 
       +\left( x \dfrac{\eta_x}{J}+y \dfrac{\eta_y}{J}+z \dfrac{\eta_z}{J}  \right)_\eta 
       +\left( x \dfrac{\zeta_x}{J}+y \dfrac{\zeta_y}{J}+z \dfrac{\zeta_z}{J}  \right)_\zeta   
     \right] .                        
\end{equation}
With this, the high-order central finite difference schemes indeed maintain free-stream preserving property. While for upwind linear and nonlinear schemes, achieving the free-stream preserving identity is not staightforward. 

\subsection{Discretization methods}
If not mentioned otherwise, the 3rd-order TVD Runge-Kutta method ~\cite{gottlieb1998total} is applied to perform the time integration. Our numerical methods are based on the conservative finite difference method, with the viscous fluxes and the convective fluxes being discretized by the 6th-order central difference scheme and the upwind high-order reconstruction schemes, respectively.
\subsubsection{Spatial discretization of the fluxes}
Without loss of generality, we consider the fluxes along the $\xi$ direction, which is indexed by $i$. In conservative FDM, let the fluxes $ \boldsymbol{\tilde{F}}_{i}$ at cell-center $i$ be the average of a primitive function $\boldsymbol{\hat{H}}(\xi)$ 
\begin{equation} \label{eq:avr}
\boldsymbol{\tilde{F}}_{i}=\dfrac{1}{\Delta \xi}\int_{i-1/2}^{i+1/2} \boldsymbol{\hat{H}}(\xi)d\xi.
\end{equation}
Then the derivative of $\boldsymbol{\tilde{F}}_{i}$ can be obtained exactly
\begin{equation} \label{eq:derivative}
\left( \dfrac{\partial \boldsymbol{\tilde{F}}}{\partial \xi} \right)_i=\dfrac{\boldsymbol{\hat{H}}(i+1/2)-\boldsymbol{\hat{H}}(i-1/2)}{\Delta \xi}
\end{equation}

The primitive function value at the cell face, say $\boldsymbol{\hat{H}}_{i+1/2}$ at $i+1/2$, can be reconstructed through the neighboring cell-center fluxes $\boldsymbol{\tilde{F}}_{i-k+1}$, $\cdots$, $\boldsymbol{\tilde{F}}_{i+k-1}$ to achieve a $(2k-1)$th-order accuracy, such as the WENO scheme, resulting in a $(2k-1)$th-order accuracy of $\partial \boldsymbol{\tilde{F}} /\partial \xi$ as well. Therefore, the derivative of the convective fluxes can be approximated by the reconstructed cell-face fluxes $\boldsymbol{\tilde{F}}_{i+1/2}$
\begin{equation} \label{eq:aderivative}
\left( \dfrac{\partial \boldsymbol{\tilde{F}}}{\partial \xi} \right)_i= \dfrac{\boldsymbol{\tilde{F}}_{i+1/2}- \boldsymbol{\tilde{F}} _{i-1/2}}{\Delta \xi}+\boldsymbol{O}(\Delta \xi^{2k-1}).
\end{equation}

\subsubsection{The characteristic-wise reconstruction scheme}
For the reconstruction of the $\boldsymbol{\tilde{F}} _{i+1/2}$, a robust way is to transform the fluxes and conservative variables into the characteristic space and apply a flux vector splitting scheme, such as the Lax-Friedrichs splitting
\begin{equation}
 \boldsymbol{ \tilde{F} }_m^{\pm}= \dfrac{1}{2} \boldsymbol{L}_{i+1/2} \cdot \left( \boldsymbol{\tilde{F}}_m \pm \boldsymbol{\lambda}_{i+1/2} \boldsymbol{\tilde{Q}}_m \right), m=i-k+1,\cdots,i+k,
\end{equation}
where $\boldsymbol{\lambda}$ is the diagonal matrix composed of the eigenvalues of the linearized Roe-average Jacobian matrix $\boldsymbol{A}_{i+1/2}=\left( \partial \boldsymbol{\tilde{F}} /\partial \boldsymbol{\tilde{Q}} \right)_{i+1/2}$. $\boldsymbol{L}_{i+1/2}$ is the left matrix composed of the corresponding eigenvectors of $\boldsymbol{A}_{i+1/2}$. The final cell-face fluxes can be computed by
 \begin{equation}
 \boldsymbol{ \tilde{F} }_{i+1/2}= \boldsymbol{R}_{i+1/2} \cdot \left( \boldsymbol{ \tilde{F} }_{i+1/2}^{+}+\boldsymbol{ \tilde{F} }_{i+1/2}^{-} \right),
\end{equation}
where $\boldsymbol{R}_{i+1/2}$ is the inverse matrix of $\boldsymbol{L}_{i+1/2}$. 

$\boldsymbol{ \tilde{F} }_{i+1/2}^{\pm}$ are reconstructed upon the specific upwind schemes, such as the 5th-order linear upwind scheme
 \begin{equation}\label{eq:linear_upwind}
 \begin{aligned}
\boldsymbol{ \tilde{F} }_{i+1/2}^{+}&=\dfrac{1}{60}\left( 2\boldsymbol{ \tilde{F} }_{i-2}^{+}-13\boldsymbol{ \tilde{F} }_{i-1}^{+}+47\boldsymbol{ \tilde{F} }_{i}^{+}+27\boldsymbol{ \tilde{F} }_{i+1}^{+}-3\boldsymbol{ \tilde{F} }_{i+2}^{+} \right) \\
\boldsymbol{ \tilde{F} }_{i+1/2}^{-}&=\dfrac{1}{60}\left( -3\boldsymbol{ \tilde{F} }_{i-1}^{-}+27\boldsymbol{ \tilde{F} }_{i}^{-}+47\boldsymbol{ \tilde{F} }_{i+1}^{-}-13\boldsymbol{ \tilde{F} }_{i+2}^{-}+2\boldsymbol{ \tilde{F} }_{i+3}^{-} \right)
\end{aligned}
\end{equation}
With respect to the nonlinear upwind reconstruction, we choose the classical 5th-order  WENO scheme~\cite{jiang_efficient_1996} to obtain the cell-face flux by
\begin{equation}\label{eq:weno}
\widetilde{f}_{i+1/2}^{\pm}=\sum\limits_{k=0}^{2}\omega_k^{\pm} q_k^{\pm},
\end{equation}
where $\widetilde{f}^{\pm}$ denotes each component of $ \boldsymbol{ \tilde{F} }^{\pm}$. Consider $\widetilde{f}^{+}$ as an example, the 3rd-order approximations for the three different sub-stencils are
\begin{equation}
\begin{aligned}
q_0^{+}&= \dfrac{1}{3}\widetilde{f}_{i-2}^{+}
          -\dfrac{7}{6}\widetilde{f}_{i-1}^{+}
          +\dfrac{7}{6}\widetilde{f}_{i}^{+},\\
q_1^{+}&=-\dfrac{1}{6}\widetilde{f}_{i-1}^{+}
          +\dfrac{5}{6}\widetilde{f}_{i}^{+}
          +\dfrac{1}{3}\widetilde{f}_{i+1}^{+},\\
q_2^{+}&= \dfrac{1}{3}\widetilde{f}_{i}^{+}
          +\dfrac{5}{6}\widetilde{f}_{i+1}^{+}
          -\dfrac{1}{6}\widetilde{f}_{i+2}^{+},               
\end{aligned}
\end{equation} 
and the corresponding nonlinear weight $\omega_k^{+}$ is proposed to be determined by
\begin{equation}\label{eq:weno_wights}
\omega_k^{+}=\dfrac{C_k}{\left( \beta_k^{+}+\epsilon \right)^n} / \sum\limits_{r=0}^{2}\dfrac{C_r}{\left( \beta_k^{+}+\epsilon \right)^n},
\end{equation} 
where $C_0=\dfrac{1}{10}$, $C_1=\dfrac{3}{5}$, $C_2=\dfrac{3}{10}$ are the optimal weights , $\epsilon=1.0\times10^{-6}$ and $n=2$. The smoothness indicators are evaluated by  
\begin{equation}
\begin{aligned}
\beta_0^{+}&= \dfrac{1}{4} \left(  \widetilde{f}_{i-2}^{+}
                                 -4\widetilde{f}_{i-1}^{+}
                                 +3\widetilde{f}_{i}^{+} \right)^2
            +\dfrac{13}{12}\left(  \widetilde{f}_{i-2}^{+}
                                 -2\widetilde{f}_{i-1}^{+}
                               +\widetilde{f}_{i}^{+} \right)^2,\\
\beta_1^{+}&= \dfrac{1}{4} \left( -\widetilde{f}_{i-1}^{+}
                               +\widetilde{f}_{i+1}^{+} \right)^2
            +\dfrac{13}{12}\left(  \widetilde{f}_{i-1}^{+}
                                 -2\widetilde{f}_{i}^{+}
                              +\widetilde{f}_{i+1}^{+} \right)^2,\\
\beta_2^{+}&= \dfrac{1}{4} \left(-3\widetilde{f}_{i}^{+}
                                 +4\widetilde{f}_{i+1}^{+}
                                -\widetilde{f}_{i+2}^{+} \right)^2
            +\dfrac{13}{12}\left(  \widetilde{f}_{i}^{+}
                                 -2\widetilde{f}_{i+1}^{+}
                              +\widetilde{f}_{i+2}^{+} \right)^2.\\                  
\end{aligned}
\end{equation} 
 
\section{Free-stream preserving strategy for upwind schemes} \label{sec:method}
In this part, we propose a novel strategy by reformulating the upwind dissipation to maintain the free-stream preserving identity for the upwind schemes. Like the existing free-stream preserving methods~\cite{nonomura_new_2015,zhu_free-stream_2019}, our method is based on the discretization of metrics with SCMM~\cite{nonomura_new_2015}, which is given by
\begin{equation}\label{metric_discret}
\begin{aligned}
\left(x_{i+1/2}\right)_{2th}&=\dfrac{1}{2}\left( x_{i}+x_{i+1} \right),\\
\left(x_{i+1/2}\right)_{4th}&=\dfrac{1}{12}\left( -x_{i-1}+7x_{i}+7x_{i+1}-x_{i+2} \right),\\
\left(x_{i+1/2}\right)_{6th}&=\dfrac{1}{60}\left( x_{i-2}-8x_{i-1}+37x_{i}+37x_{i+1}-8x_{i+2}+x_{i+3} \right),\\
\left(x_{i+1/2}\right)_{8th}&=\dfrac{1}{840}\left( -3x_{i-3}+29x_{i-2}-139x_{i-1}+533x_{i}+533x_{i+1}-139x_{i+2}+29x_{i+3}-3x_{i+4} \right),\\
\left( \dfrac{\partial x}{\partial \xi} \right)_i&=x_{i+1/2}-x_{i-1/2}.
\end{aligned}
\end{equation} 

In the following discussions, without loss of generality, the 5th-order linear upwind and WENO schemes are considered to reconstruct the cell-face fluxes with this suggested free-stream preserving strategy. 

\subsection{Linear upwind scheme}
As derived in Ref.~\cite{zhu_free-stream_2019}, the 5th-order linear upwind scheme in Eq.~\eqref{eq:linear_upwind} can be formulated by a 6th-order central part where only the cell-center fluxes are included, and the 5th-order dissipation term which only consists of $\boldsymbol{\tilde{Q}}$,
\begin{equation}\label{Linear_upwind_standard}
\begin{aligned}
 \boldsymbol{ \tilde{F} }_{i+1/2}=&\dfrac{1}{60} \left(\boldsymbol{ \tilde{F} }_{i-2}-8\boldsymbol{ \tilde{F} }_{i-1}+37\boldsymbol{ \tilde{F} }_{i}+37\boldsymbol{ \tilde{F} }_{i+1}-8\boldsymbol{ \tilde{F} }_{i+2}+\boldsymbol{ \tilde{F} }_{i+3} \right)\\
 &+\dfrac{1}{60}\sum\limits_s \boldsymbol{R}_{i+1/2}^s  \lambda^s \boldsymbol{L}_{i+1/2}^s \left( \boldsymbol{\tilde{Q}}_{i-2}-5\boldsymbol{\tilde{Q}}_{i-1}+10\boldsymbol{\tilde{Q}}_{i}-10\boldsymbol{\tilde{Q}}_{i+1}+5\boldsymbol{\tilde{Q}}_{i+2}-\boldsymbol{\tilde{Q}}_{i+3} \right).
\end{aligned}
\end{equation}
Therefore, if the 6th-order central scheme in Eq.~\eqref{metric_discret} is employed to discretize the metrics and Jacobian, only the local dissipation may violate the free-stream preserving property. Such issue can be addressed by reformulating the upwind dissipation portion of each $\boldsymbol{\tilde{F}}_m^{\pm}$ during the reconstruction,
\begin{equation}\label{eq:Linear_upwind_LFspliting}
\begin{aligned}
\boldsymbol{\tilde{F}}_m^{\pm}&=\dfrac{1}{2} \boldsymbol{L}_{i+1/2}\left[ \boldsymbol{\tilde{F}}_m \pm \boldsymbol{\lambda} \left(\dfrac{\boldsymbol{Q}}{J}\right)_m \right], m=i-2,\cdots,i+3\\
&\approx \dfrac{1}{2} \boldsymbol{L}_{i+1/2} \left\lbrace \boldsymbol{\tilde{F}}_m \pm \boldsymbol{\lambda} \left[ \left(\dfrac{\boldsymbol{Q}-\boldsymbol{Q}^*}{J}\right)_m +\boldsymbol{Q}^*\left(\dfrac{1}{J}\right)_{i+1/2} \right] \right\rbrace,
\end{aligned}
\end{equation}
where $\boldsymbol{Q^*}$ is the reference cell-face conservative variables along the stencil, proposed by
\begin{equation}\label{Q_star}
 \boldsymbol{Q^*}=\boldsymbol{\tilde{Q}}_{i+1/2}/\left(\dfrac{1}{J}\right)_{i+1/2},
\end{equation} 
and
\begin{equation}\label{QdJ_star}
 \boldsymbol{\tilde{Q}}_{i+1/2}=\dfrac{1}{60}\left( \boldsymbol{\tilde{Q}}_{i-2}-8\boldsymbol{\tilde{Q}}_{i-1}+37\boldsymbol{\tilde{Q}}_{i}+37\boldsymbol{\tilde{Q}}_{i+1}-8\boldsymbol{\tilde{Q}}_{i+2}+\boldsymbol{\tilde{Q}}_{i+3} \right),
\end{equation}
\begin{equation}\label{1dJ}
\left(\dfrac{1}{J}\right)_{i+1/2}=\dfrac{1}{60} \left[ \left(\dfrac{1}{J}\right)_{i-2}-8\left(\dfrac{1}{J}\right)_{i-1}+37\left(\dfrac{1}{J}\right)_{i}+37\left(\dfrac{1}{J}\right)_{i+1}-8\left(\dfrac{1}{J}\right)_{i+2}+\left(\dfrac{1}{J}\right)_{i+3} \right].
\end{equation} 
Substituting Eq.~\eqref{eq:Linear_upwind_LFspliting} into Eq.~\eqref{eq:linear_upwind} yields
\begin{equation} \label{eq:final_form_linear_upwind}
\begin{aligned} 
 \boldsymbol{ \tilde{F} }_{i+1/2}^{'}=&\dfrac{1}{60} \left(\boldsymbol{ \tilde{F} }_{i-2}-8\boldsymbol{ \tilde{F} }_{i-1}+37\boldsymbol{ \tilde{F} }_{i}+37\boldsymbol{ \tilde{F} }_{i+1}-8\boldsymbol{ \tilde{F} }_{i+2}+\boldsymbol{ \tilde{F} }_{i+3} \right)\\
 &+\dfrac{1}{60}\sum\limits_s \boldsymbol{R}_{i+1/2}^s  \lambda^s \boldsymbol{L}_{i+1/2}^s \left[ \left(\boldsymbol{\tilde{Q}}_{i-2}-\boldsymbol{\tilde{Q}}_{i-2}^*\right)-5\left(\boldsymbol{\tilde{Q}}_{i-1}-\boldsymbol{\tilde{Q}}_{i-1}^*\right)+10\left(\boldsymbol{\tilde{Q}}_{i}-\boldsymbol{\tilde{Q}}_{i}^*\right) \right.\\ 
 &  \left. -10\left(\boldsymbol{\tilde{Q}}_{i+1}-\boldsymbol{\tilde{Q}}_{i+1}^*\right)+5\left(\boldsymbol{\tilde{Q}}_{i+2}-\boldsymbol{\tilde{Q}}_{i+2}^*\right)-\left(\boldsymbol{\tilde{Q}}_{i+3}-\boldsymbol{\tilde{Q}}_{i+3}^*\right) \right]\\
=& \boldsymbol{ \tilde{F} }_{1+1/2}-\dfrac{1}{60} \Delta^5 \left( \dfrac{1}{J} \right)    \sum\limits_s \boldsymbol{R}_{i+1/2}^s  \lambda^s \boldsymbol{L}_{i+1/2}^s \boldsymbol{Q}^*,
\end{aligned}
\end{equation}
where 
\begin{equation}\label{eq:linear_upwind_adding}
\begin{aligned}
-\Delta^5 \left( \dfrac{1}{J} \right)  &=-\left[ \left(\dfrac{1}{J}\right)_{i-2}-5\left(\dfrac{1}{J}\right)_{i-1}+10\left(\dfrac{1}{J}\right)_{i}-10\left(\dfrac{1}{J}\right)_{i+1}+5\left(\dfrac{1}{J}\right)_{i+2}-\left(\dfrac{1}{J}\right)_{i+3} \right] \\
&=\left( \dfrac{1}{J^6}\dfrac{\partial^5 J}{\partial \xi^5} \right)_{i+1/2}  \Delta \xi ^5+O\left(\Delta \xi ^6 \right).
\end{aligned}
\end{equation}

Upon this modification, the reformulated dissipation of the proposed linear upwind scheme can maintain the free-stream preserving identity essentially as $\boldsymbol{Q}-\boldsymbol{Q}^*$ vanish strictly under the free-stream condition. In addition, as shown in Eq.~\eqref{eq:Linear_upwind_LFspliting}, the present method recover the original linear upwind scheme Eq. \eqref{Linear_upwind_standard} on a uniform Cartesian grid. Essentially, the cell-face fluxes $\boldsymbol{\tilde{F}}_{i+1/2}^{'}$ obtained by the proposed strategy is an approximation of the standard 5th-order linear upwind scheme by adding an extra 5th-order dissipation term, as elaborated in Eq.~\eqref{eq:final_form_linear_upwind}.

\subsection{WENO scheme}
Inspired by the above linear upwind scheme, we split the cell-center fluxes of the WENO reconstruction in the local stencil into two parts as 
\begin{equation} \label{flux_split_detail}
\begin{aligned}
 \boldsymbol{\tilde{F}}_m = &\left( \boldsymbol{\tilde{F}}_m-\boldsymbol{\tilde{F}}_m^{*} \right)+\boldsymbol{\tilde{F}}_{m}^{*} , \quad m=i-2, \cdots, i+3 \\
 =& \left( \boldsymbol{F}_m-\boldsymbol{F^*} \right) \left(\dfrac{\xi_x}{J}\right)_m 
    +\boldsymbol{F^*} \left(\dfrac{\xi_x}{J} \right)_{m} \\
  &+\left( \boldsymbol{G}_m-\boldsymbol{G^*} \right) \left(\dfrac{\xi_y}{J}\right)_m 
    +\boldsymbol{G^*} \left(\dfrac{\xi_y}{J} \right)_{m} \\
  &+\left( \boldsymbol{H}_m-\boldsymbol{H^*} \right) \left(\dfrac{\xi_z}{J}\right)_m 
    +\boldsymbol{H^*} \left(\dfrac{\xi_z}{J} \right)_{m}.
\end{aligned}
\end{equation} 
$\boldsymbol{F^*}$, $\boldsymbol{G^*}$ and $\boldsymbol{H^*}$ are the corresponding fluxes of the reference cell-face conservative variables $\boldsymbol{Q^*}$ in Eq.~\eqref{Q_star}, and they can be considered as the approximations of $\boldsymbol{F}$, $\boldsymbol{G}$ and $\boldsymbol{H}$ at cell-face $i+1/2$. Here, the reconstruction of $\boldsymbol{\tilde{F}}_m-\boldsymbol{\tilde{F}}_m^{*}$ always maintains free-stream preserving identity regardless of the particular formulation of the employed WENO schemes. For the reconstruction of $\boldsymbol{\tilde{F}}_{m}^{*}$, the central scheme can achieve this property, as it retains the sufficient condition of Deng and Abe ~\cite{ DENG201390,abe2014geometric}. Thus, we achieve
\begin{equation} \label{eq:my_point1}
\begin{aligned}
 \boldsymbol{\tilde{F}}_{i+1/2}^{'}&=WENO \left( \boldsymbol{\tilde{F}}_{i-2}-\boldsymbol{\tilde{F}}_{i-2}^*, \cdots, \boldsymbol{\tilde{F}}_{i+3}-\boldsymbol{\tilde{F}}_{i+3}^* \right) +CENTRAL \left( \boldsymbol{\tilde{F}}_{i-2}^*, \cdots,\boldsymbol{\tilde{F}}_{i+3}^*\right) \\
                                             &=WENO \left( \boldsymbol{\tilde{F}}_{i-2}-\boldsymbol{\tilde{F}}_{i-2}^*, \cdots, \boldsymbol{\tilde{F}}_{i+3}-\boldsymbol{\tilde{F}}_{i+3}^* \right)+ \boldsymbol{\tilde{F}}_{i+1/2}^*,
\end{aligned}
\end{equation} 
where $WENO(*)$ operation stands for the reconstruction in Eq.~\eqref{eq:weno}, and $CENTRAL(*)$ is the central part of Eq.~\eqref{Linear_upwind_standard} which leads to
\begin{equation} 
\begin{aligned}
 \boldsymbol{\tilde{F}}_{i+1/2}^{*}=  \boldsymbol{F^*}  \left(\dfrac{\xi_x}{J}\right)_{i+1/2}  
	+\boldsymbol{G^*}  \left(\dfrac{\xi_y}{J}\right)_{i+1/2} 
	+\boldsymbol{H^*}  \left(\dfrac{\xi_z}{J}\right)_{i+1/2}.
\end{aligned}
\end{equation} 
 Like Eq.~\eqref{QdJ_star} and Eq.~\eqref{1dJ}, the cell-face metrics, say $ \left(\xi_x/J\right)_{i+1/2}$, can be approximated by the 6th-order scheme as
\begin{equation}\label{xidJ}
\left(\dfrac{\xi_x}{J}\right)_{i+1/2}=\dfrac{1}{60} \left[ \left(\dfrac{\xi_x}{J}\right)_{i-2}-8\left(\dfrac{\xi_x}{J}\right)_{i-1}+37\left(\dfrac{\xi_x}{J}\right)_{i}+37\left(\dfrac{\xi_x}{J}\right)_{i+1}-8\left(\dfrac{\xi_x}{J}\right)_{i+2}+\left(\dfrac{\xi_x}{J}\right)_{i+3} \right].
\end{equation} 
 
For simplicity, as the $CENTRAL$ part $\boldsymbol{\tilde{F}}_{i+1/2}^*$ is a constant vector, it can be moved into the $WENO$ part without changing the final reconstructed fluxes $\boldsymbol{\tilde{F}}_{i+1/2}^{'}$. Therefore, we reformulate Eq.~\eqref{eq:my_point1} by
\begin{equation} \label{eq:my_point2}
\begin{aligned}
 \boldsymbol{\tilde{F}}_{i+1/2}^{'}&=WENO \left( \boldsymbol{\tilde{F}}_{i-2}-\boldsymbol{\tilde{F}}_{i-2}^*+\boldsymbol{\tilde{F}}_{i+1/2}^*, \cdots, \boldsymbol{\tilde{F}}_{i+3}-\boldsymbol{\tilde{F}}_{i+3}^*+ \boldsymbol{\tilde{F}}_{i+1/2}^* \right).
\end{aligned}
\end{equation}
In addtion, from another perspective, this treatment can be considered as an approximation of the cell-center fluxes before the WENO reconstruction,
\begin{equation} \label{eq:my_point3}
\begin{aligned}
 \boldsymbol{\tilde{F}}_m \approx &\boldsymbol{\tilde{F}}_{m}-\boldsymbol{\tilde{F}}_{m}^*+\boldsymbol{\tilde{F}}_{i+1/2}^*\\
 = & \left( \boldsymbol{f}_m-\boldsymbol{f^*} \right) \left(\dfrac{\xi_x}{J}\right)_m 
                                                                                     +\boldsymbol{f^*} \left(\dfrac{\xi_x}{J} \right)_{i+1/2} \\
														&+\left( \boldsymbol{g}_m-\boldsymbol{g^*} \right) \left(\dfrac{\xi_y}{J}\right)_m 
                                                                                     +\boldsymbol{g^*} \left(\dfrac{\xi_y}{J} \right)_{i+1/2} \\
														&+\left( \boldsymbol{h}_m-\boldsymbol{h^*} \right) \left(\dfrac{\xi_z}{J}\right)_m 
                                                                                     +\boldsymbol{h^*} \left(\dfrac{\xi_z}{J} \right)_{i+1/2}.
\end{aligned}
\end{equation}

\begin{remark}\label{remark1}
For a given the uniform Cartesian grid, the present method recovers to the standard WENO scheme, as $\boldsymbol{\tilde{F}}_{m}^* \left(m=i-2,\cdots,i+3\right)$ equal to $\boldsymbol{\tilde{F}}_{i+1/2}^*$.
\end{remark}

In addtion, as deailed in Appendix, the reconstructed cell-face fluxes in Eq. \eqref{eq:my_point2} can be expressed as the combination of the central and dissipation parts. 
We note that this approximation only reformulates the upwind dissipation, without altering the central part of the standard WENO schemes, as $\boldsymbol{\hat{F}}_{i+1/2}$ in Eq.~\eqref{eq:WENO_hatF} of the Appendix equals zero if the metrics and Jacobian are discretized with the 6th-order central scheme.  

\begin{remark}\label{remark2}
Unlike previous FP methods in e.g. Refs. \cite{nonomura_new_2015,zhu_free-stream_2019}, this method does not change the standard form of the WENO scheme, as in Eq.~\eqref{eq:my_point2}.
In this way, the extension to other WENO schemes, such as the standard 7th-order WENO scheme~\cite{shu1998essentially} and the WENOCU6 scheme~\cite{HU20108952}, is straightforward.
\end{remark}

\begin{remark}\label{remark3}
The numerical strategy in Eq.~\eqref{eq:my_point1} and Eq.~\eqref{eq:my_point2} removes the constraint that a consistent discretization is applied for the metrics and fluxes, 
as long as the conservative quantities in Eq.~\eqref{QdJ_star} use the same scheme with the that of metrics.
\end{remark}

\section{Numerical tests on curvilinear grids} \label{sec:testcases}
Several problems, including the free-stream, isotropic vortex convection, double Mach reflection, subsonic flow past the 30P30N three-element airfoil and transonic flow around the three-dimensinal (3D) ONERA M6 wing etc. are conducted to check the performance of the proposed free-stream preserving method on the curvilinear grids. The global Lax-Friedrichs flux splitting is employed to the double Mach reflection and supersonic flow pass a cylinder while other verifications are performed with the local one. In addition, the geometric metrics are obtained by the 2nd- to 8th-order central schemes. In the followings, WENO5 and WENO7 are the standard 5th-order and 7th-order WENO scheme~\cite{jiang_efficient_1996,shu1998essentially}. WENOCU6 denotes the 6th-order WENOCU6 scheme~\cite{HU20108952}. Besides, the WENO schemes with the proposed free-stream preserving strategy are denoted by WENOX-MY, which presents the Xth-order WENO together with Yth-order metrics and Jacobian, respectively. For example, WENO7-M6 means the 7th-order WENO and 6th-order metrics and Jacobian, and Linear-Upwind5-M6 denotes the suggested 5th-order linear upwind scheme coupled with the 6th-order metrics and Jacobian. 

\subsection{Simple test cases}
\subsubsection{Free-streaming problem}\label{sec:Free-stream}
First, we consider the wavy and random grids to test the free-stream preservation, as shown in Fig.~\ref{wavy_grid}. The wavy case is defined in the domain $(x,y)\in [-10,10]\times[-10,10]$ by
\begin{equation}\label{eq:wavy_grid}
\begin{aligned}
x_{i,j}=x_{min}+\Delta{x_0} \left[ (i-1)+A_x sin \left( \dfrac{n_{xy}\pi(j-1)\Delta{y_0}}{L_y} \right)\right] \\
y_{i,j}=y_{min}+\Delta{y_0} \left[ (j-1)+A_y sin \left( \dfrac{n_{yx}\pi(i-1)\Delta{x_0}}{L_x} \right)\right],
\end{aligned}
\end{equation}
where $L_x=L_y=20$, $x_{min}=-L_x/2$, $y_{min}=-L_y/2$, $A_x\Delta{x}=0.6$, $A_y\Delta{y}=0.6$, and $n_{xy}=n_{yx}=8$. And the grid points of the random case is generated in a random direction with $20\%$ of the original Cartesian grid size. A coarse grid resolution of $21\times21$ is chosen here to highlight the differences between the standard WENO and the proposed free-stream preserving WENO schemes, for both the wavy and random grids.
\begin{figure} [htbp]
    \center 
    \includegraphics[width=0.8\textwidth]{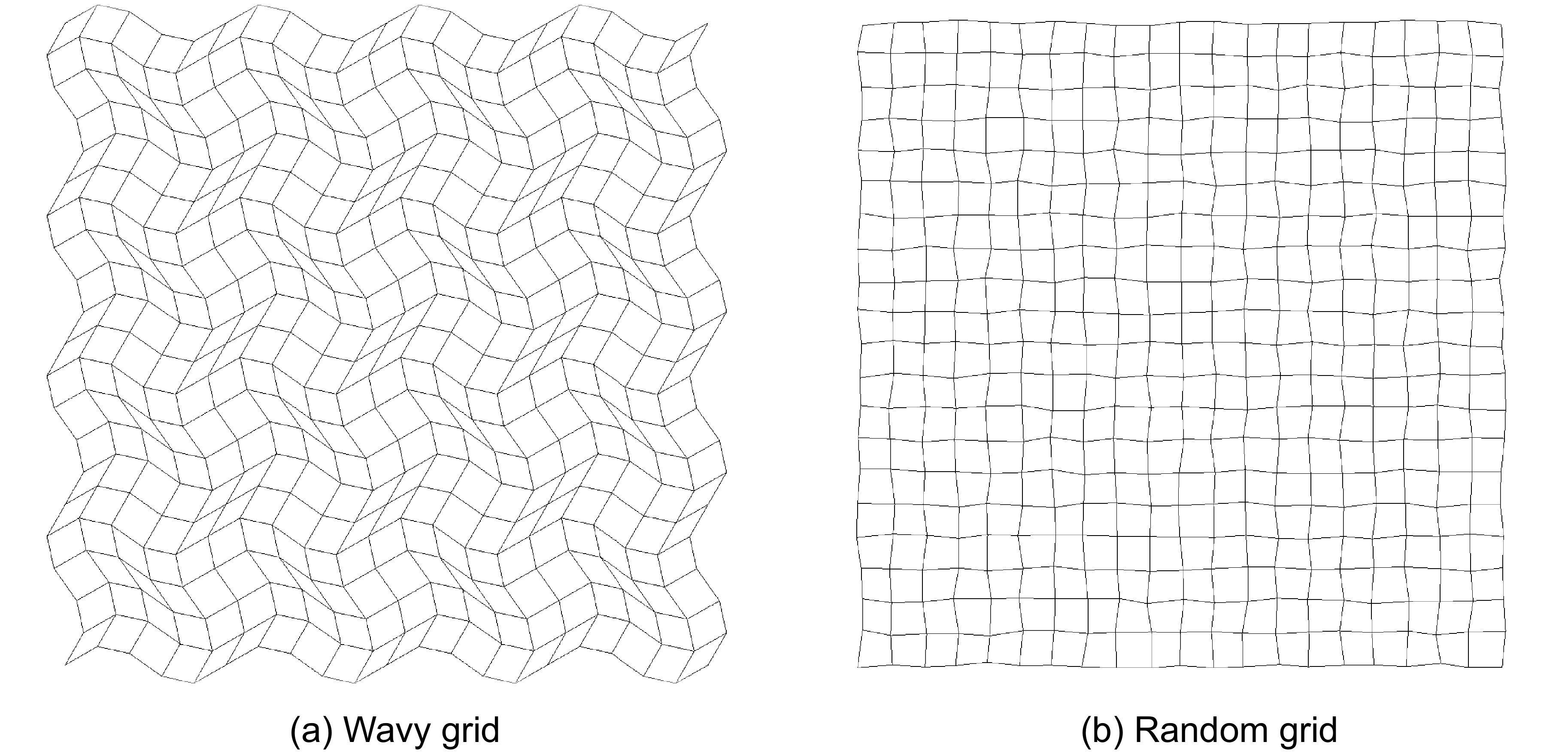}
    \caption{The wavy and random grids at a resolution of $21 \times 21$. }
    \label{grid41}
\end{figure}

The initial condition is uniform with a $Ma=0.5$ flow in $x$ direction,
\begin{equation}
u=0.5,v=0,p=1,\rho=\gamma,
\end{equation}
where $\gamma=1.4$ is the specific heat ratio. The time step is $0.2$, and the average ($L_2$) and maximal ($L_\infty$) errors of the velocity components $v$ for the two grids are estimated at $t=20$. As show in Table~\ref{freestream_L2}, both the $L_2$ and $L_\infty$ errors of the original WENO5 and WENO7 schemes have been significantly reduced to almost the machine zero by our proposed FP method (WENO5-M6, WENO7-M6, WENO7-M8), indicating that our method indeed preserve the free-streaming identity.

\begin{table}[htbp]	
   \centering
	\caption{The $L_2$ and $L_{\infty}$ errors of the $v$ component on the wavy and random grids at a resolution of $21\times21$.}	 
	
	\begin{tabular}{ccccc}
   \hline
    \multirow{2}*{Method}  & \multicolumn{2}{c}{Wavy grid}   &\multicolumn{2}{c}{Random grid} \\
   \cline{2-5}          
    ~ &$L_2$ error  &$L_{\infty}$ error &$L_2$ error  &$L_{\infty}$ error  \\ 
   \hline 
   WENO5         &$2.45\times10^{-2}  $  &$4.72\times10^{-2}$  &$1.29\times10^{-2}$    &$4.41\times10^{-2}$ \\  
   WENO7         &$1.03\times10^{-2}  $  &$1.98\times10^{-2}$  &$1.57\times10^{-2}$    &$5.08\times10^{-2}$ \\   
   WENO5-M6   &$5.58\times10^{-16}$ &$2.05\times10^{-15}$ &$7.60\times10^{-16}$  &$2.70\times10^{-15}$ \\
   WENO7-M6   &$5.88\times10^{-16}$ &$1.92\times10^{-15}$ &$8.35\times10^{-16}$  &$2.54\times10^{-15}$ \\
   WENO7-M8   &$4.90\times10^{-16}$ &$1.94\times10^{-15}$ &$6.20\times10^{-16}$  &$1.96\times10^{-15}$ \\
		\hline
	\end{tabular}
	\label{freestream_L2}
\end{table}
\subsubsection{Isotropic vortex}\label{sec:Isotropic vortex}
The two-dimensional ($2D$) moving isotropic vortex problem on wavy and random grids is considered to evaluate the vortex preservation property. In this case, the fluid is treated as an ideal gas with the specific heat ratio $\gamma=1.4$. And the initial condition is the superposition of an isotropic vortex (centered at $(x_c,y_c)=(0,0)$) and an uniform flow of Mach $0.5$, with the perturbations of the velocity, temperature and entropy being
\begin{equation}
\begin{aligned}
&(\delta u,\delta v)=\epsilon \tau e^{\alpha(1-\tau^2)} (sin{\theta},-cos{\theta}), \\
&\delta{T}=-\dfrac{(\gamma-1)\epsilon^2}{4\alpha\gamma}e^{2\alpha(1-\tau^2)}, \\
&\delta{S}=0,
\end{aligned}
\end{equation}
where $\alpha=0.204$, $\tau=r/r_c$, and $r=\sqrt{(x-x_c)^2+(y-y_c)^2}$. Here $r_c=1.0$ and $\epsilon=0.02$ denote the vortex core length and strength, respectively. $T=p/\rho$ is the temperature, and $S=p/\rho^\gamma$ is the entropy. The periodic boundary condition is imposed, and the vortex moves back to the original position at $t=40$.

\begin{figure} [htbp]
    \center 
    \includegraphics[width=0.75\textwidth]{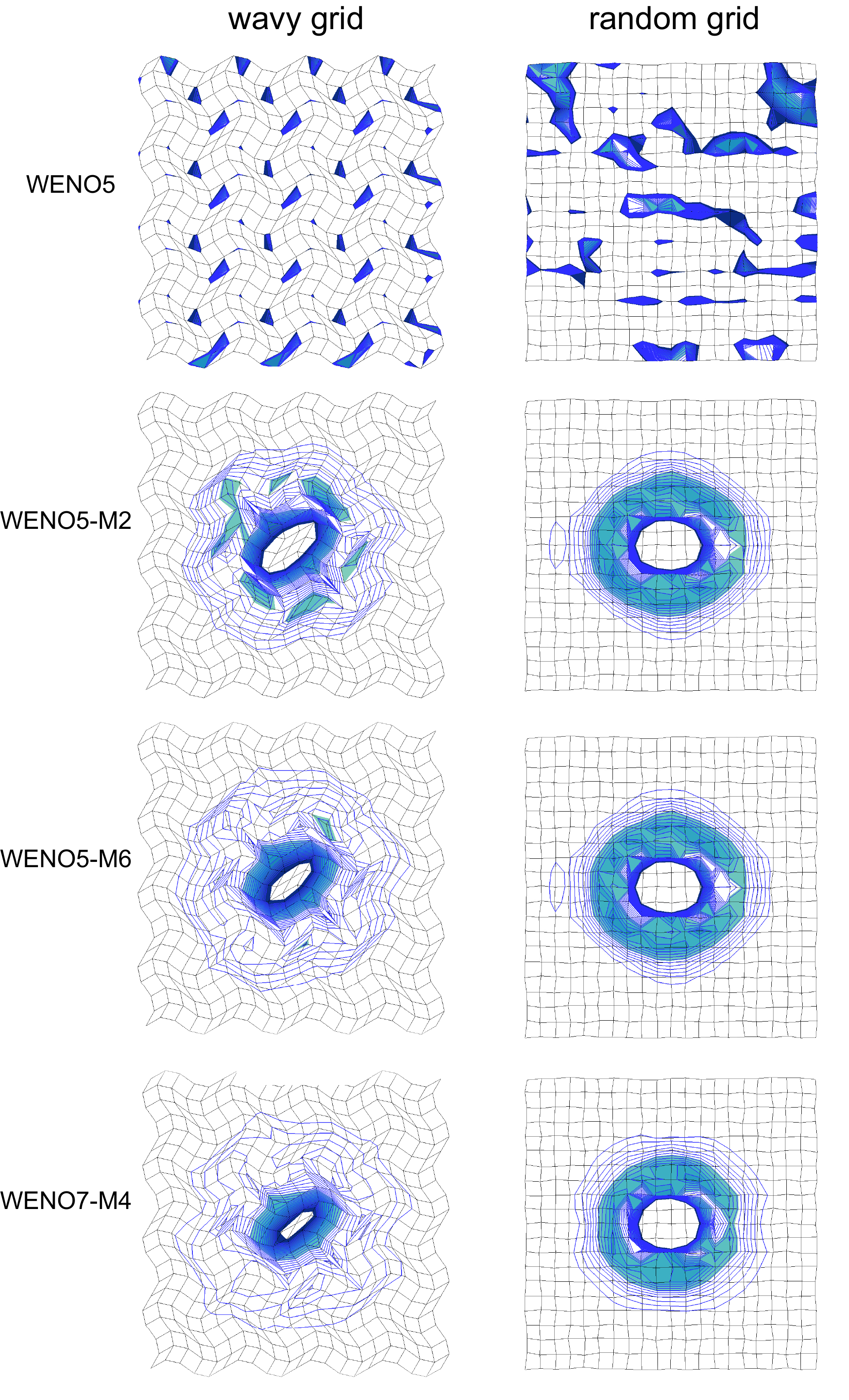}
    \caption{The vorticity magnitude contours ranging from $0.0$ to $0.006$ of the 2D moving vortex on both the wavy and random grids at $t=40$. }
    \label{vortex_grid}
\end{figure}

As shown in Fig. \ref{vortex_grid}, the vorticity magnitude contours of the original WENO5 scheme do not maintain the vortical structure on a relatively coarse resolution ($21\times21$) for both the wavy and random grids in Fig. \ref{grid41}, due to the free-stream violation errors. While this vortex clearly exhibits in the numerical results of the proposed free-stream schemes. Next, we assess the grid convergence rate of the proposed schemes on the wavy grids with the resolution of $21\times21$, $41\times41$, $81\times81$, $161\times161$ and $321\times321$. And the time step $\Delta{t}$ of each grid is decreased until the $L_2$ and $L_{\infty}$ errors are invariant to eliminate the time integration errors, as suggested by Refs. \cite{nonomura_new_2015} and \cite{zhu_free-stream_2019}. The $L_2$ and $L_{\infty}$ errors of $v$ on those wavy grids are listed in Table~\ref{wavy_grid}, as well as the corresponding convergence rates. The numerical results indicate that the schemes with consistent flux reconstruction and metric reconstruction, i.e. the order of the central part of the flux reconstruction method is the same with the order of the metric reconstruction method, including the Linear-Upwind5-M6, WENO5-M6, and WENO7-M8, show their corresponding formal order of accuracy with the resolution increasing. On the other hand, for WENO5-M2, WENO5-M4 and WENO7-M4, when the resolution is increased, the obtained order of accuracy is dominated by the metric reconstruction portion whose order is lower than that of the the flux reconstruction. However, we note here that for coarse and moderate grid resolutions, from $21 \times 21$ to $81 \times 81$, these schemes have approximately the same accuracy with that of WENO5-M6 and WENO7-M8. Thus, we conclude that the discretization errors of the metrics do not have a significant effect on the accuracy of the under-resolved simulations. 

\begin{table}[htbp]	
   \centering
	\caption{The $L_2$ and $L_{\infty}$ errors of the $v$ component and the corresponding convergence rates on the wavy grid.}
\resizebox{0.95\textwidth}{!}{	 	
	\begin{tabular}{cccccc}
	    \hline
    Method &Grid size &$L_2$ error  &Convergence rate  &$L_{\infty}$ error  &Convergence rate \\ 	
		\hline 	
		Linear-Upwind5-M6&$21 \times 21$    &$1.83\times10^{-3}$ &-       &$1.32\times10^{-2}$  &- \\
		          &$41\times 41$     &$3.53\times10^{-4}$ &2.37    &$3.32\times10^{-3}$  &1.99 \\
		          &$81\times 81$     &$1.69\times10^{-5}$ &4.38    &$1.60\times10^{-4}$  &4.38 \\
		          &$161\times 161$   &$8.77\times10^{-7}$ &4.27    &$1.08\times10^{-5}$  &3.89 \\
	  	          &$321\times 321$   &$3.02\times10^{-8}$ &4.86    &$4.06\times10^{-7}$  &4.73 \\	         
	  	\hline  	
		WENO5-M2&$21 \times 21$    &$2.12\times10^{-3}$ &-       &$1.53\times10^{-2}$  &- \\
		          &$41\times 41$     &$4.75\times10^{-4}$ &2.16    &$4.28\times10^{-3}$  &1.84 \\
		          &$81\times 81$     &$1.63\times10^{-5}$ &4.86    &$1.22\times10^{-4}$  &5.13 \\
		          &$161\times 161$   &$4.29\times10^{-6}$ &1.93    &$4.13\times10^{-5}$  &1.56 \\
	  	          &$321\times 321$   &$1.12\times10^{-6}$ &1.94    &$1.01\times10^{-5}$  &2.03 \\                    
	  	\hline  	
		WENO5-M4&$21 \times 21$    &$2.35\times10^{-3}$ &-       &$1.66\times10^{-2}$  &- \\
		          &$41\times 41$     &$5.40\times10^{-4}$ &2.12    &$4.94\times10^{-3}$  &1.75 \\
		          &$81\times 81$     &$1.72\times10^{-5}$ &4.97    &$1.65\times10^{-4}$  &4.90 \\
		          &$161\times 161$   &$8.71\times10^{-7}$ &4.30    &$1.07\times10^{-5}$  &3.95 \\
	  	          &$321\times 321$   &$2.99\times10^{-8}$ &4.86    &$3.98\times10^{-7}$  &4.75 \\             
		\hline  	
		WENO5-M6&$21 \times 21$    &$2.42\times10^{-3}$ &-       &$1.70\times10^{-2}$  &- \\
		          &$41\times 41$     &$5.47\times10^{-4}$ &2.15    &$5.02\times10^{-3}$  &1.75 \\
		          &$81\times 81$     &$1.74\times10^{-5}$ &4.97    &$1.68\times10^{-4}$  &4.90 \\
		          &$161\times 161$   &$8.77\times10^{-7}$ &4.31    &$1.08\times10^{-5}$  &3.96 \\
	  	          &$321\times 321$   &$3.02\times10^{-8}$ &4.86    &$4.06\times10^{-7}$  &4.73 \\	  	          
		\hline  	
		WENO7-M4&$21 \times 21$    &$2.59\times10^{-3}$ &-       &$1.75\times10^{-2}$  &- \\
		          &$41\times 41$     &$6.42\times10^{-4}$ &2.01    &$5.64\times10^{-3}$  &1.63 \\
		          &$81\times 81$     &$6.67\times10^{-6}$ &6.59    &$7.01\times10^{-5}$  &6.33 \\
		          &$161\times 161$   &$5.34\times10^{-8}$ &6.96    &$8.35\times10^{-7}$  &6.39 \\
                &$321\times 321$   &$1.43\times10^{-9}$ &5.22   &$1.64\times10^{-8}$  &5.67 \\                
       \hline  	
		WENO7-M8&$21 \times 21$      &$2.70\times10^{-3}$   &-          &$1.80\times10^{-2}$  &- \\
		                    &$41\times 41$       &$6.55\times10^{-4}$   &2.04    &$5.73\times10^{-3}$  &1.65 \\
		           		  &$81\times 81$       &$6.82\times10^{-6}$   &6.59    &$7.24\times10^{-5}$  &6.31 \\
		          			  &$161\times 161$   &$5.19\times10^{-8}$   &7.04    &$7.59\times10^{-7}$  &6.58 \\
                			  &$321\times 321$   &$4.49\times10^{-10}$ &6.85    &$6.89\times10^{-9}$  &6.78 \\
		\hline  		
	\end{tabular}
	\label{wavy_grid}
	}
\end{table}

\subsection{Inviscid wall-bounded flows}
\subsubsection{Double Mach reflection}\label{sec:DM}
\begin{figure} 
    \center 
    \includegraphics[width=\textwidth]{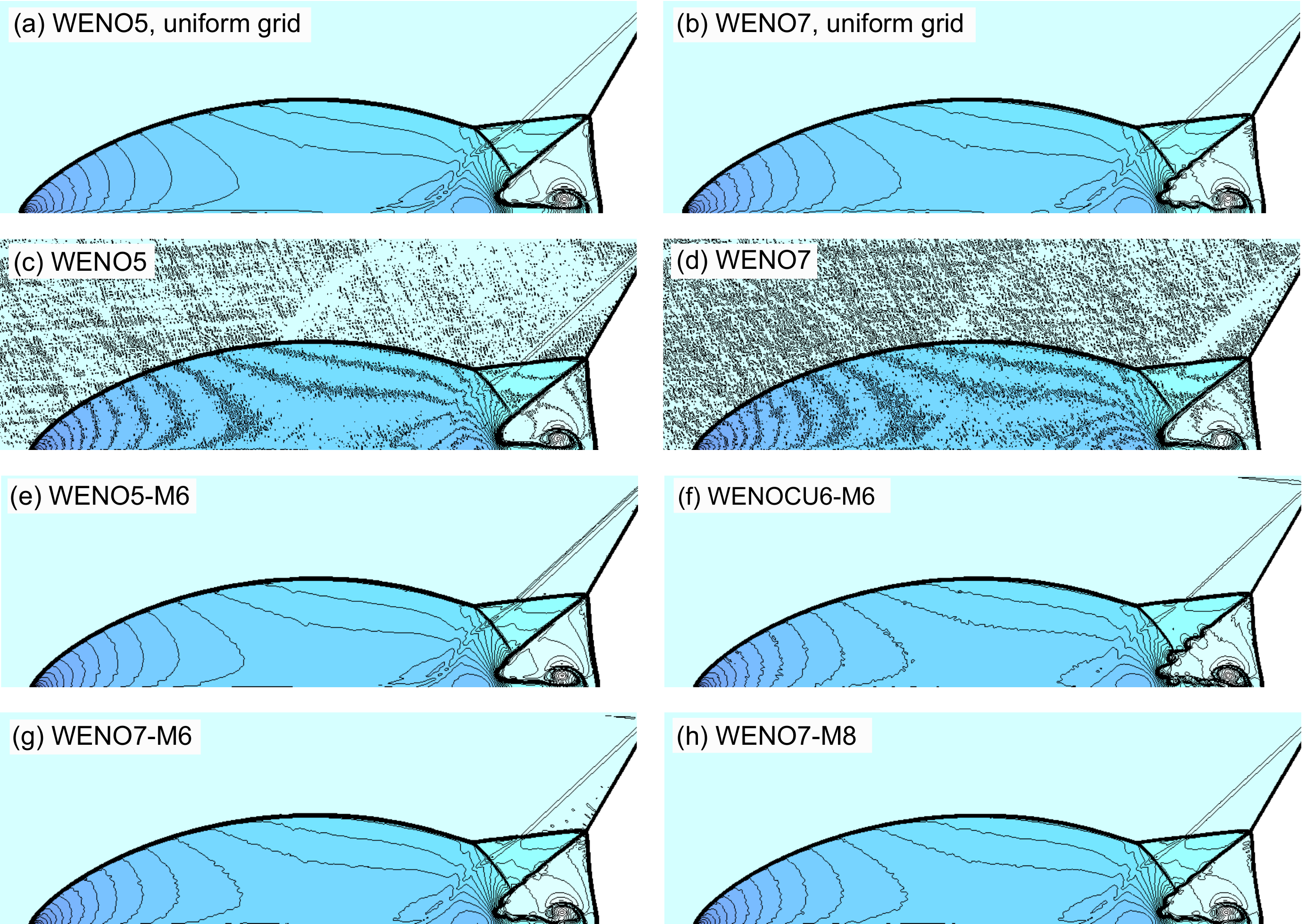}
    \caption{The density contours ranging from $1.25$ to $21.5$ of the double Mach reflection problem with the grid ($961\times241$ cells) disturbed by $5\%$ randomness.}
    \label{double_mach1}
\end{figure}

\begin{figure} 
    \center 
    \includegraphics[width=\textwidth]{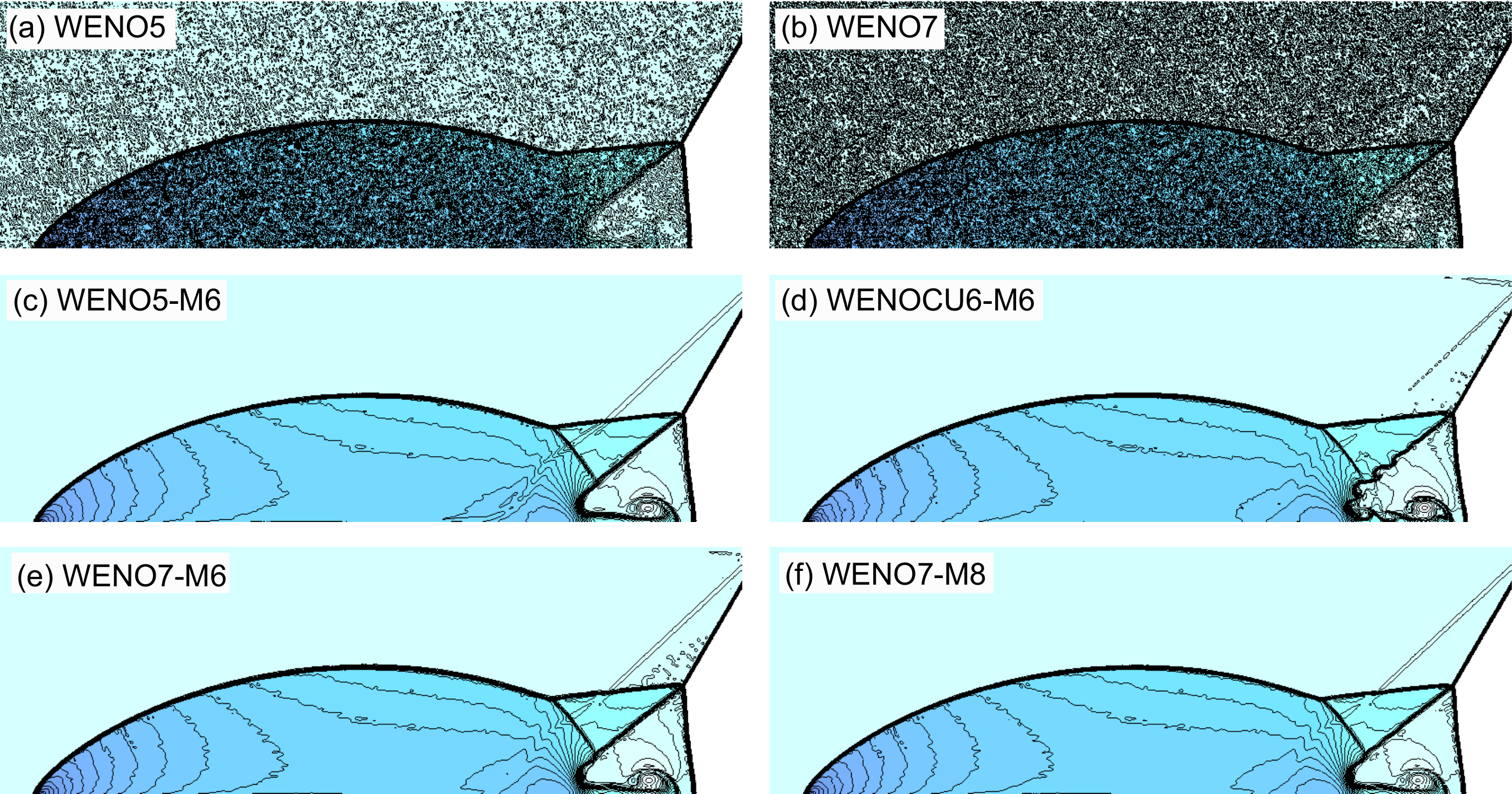}
    \caption{The density contours ranging from $1.25$ to $21.5$ of the double Mach reflection problem with the grid ($961\times241$ cells) disturbed by $20\%$ randomness.}
    \label{double_mach2}
\end{figure}

Now we test the shock-capturing capbility of the proposed numerical methods by the double Mach problem~\cite{woodward_numerical_1984} with the following initial condition
\begin{equation}
\begin{aligned}
 \left( \rho,u,v,p \right)= 
 \begin{cases}
 (1.4,0,0,1.0)                & \text{if} ~ x-y \tan(\pi/6) \geq 1/6 ,  \\
 (8.0, 7.1447, -4.125, 116.5) & \text{if} ~ x-y \tan(\pi/6) < 1/6  .
 \end{cases}
\end{aligned}
\end{equation}
The grids resolution is $961\times241$, with the randomness of the grid points in Fig. \ref{double_mach1} and Fig. \ref{double_mach2} being $5\%$ and $20\%$, respectively. The computation is performed until $t=0.2$ with a CFL number of 0.6.

As shown in Fig. \ref{double_mach1} (a-d), although the standard WENO5 and WENO7 schemes perform well on the uniform Cartesian grid, they exhibit large errors (as indicated by the spurious oscillations) in the grid with $5\%$ randomness, even in the smooth regions. While the WENO5-M6, WENOCU6-M6, WENO7-M6, and WENO7-M8 schemes achieve smoother and more continuous flow fields in Fig. \ref{double_mach1} (e-h), indicating the propose FP methods is able to accurately and robustly capture the shock and vortical structures in the double Mach reflection problem on such highly disturbed grid.
We also notice that increasing the randomness from $5\%$ to $20\%$ does not affect the numerical results in the FP methods, although larger errors are observed in the original schemes, as shown in Fig. \ref{double_mach2}.


\subsubsection{Supersonic flow past a cylinder} \label{sec:sphere}
\begin{figure} [htbp]
    \center 
    \includegraphics[width=\textwidth]{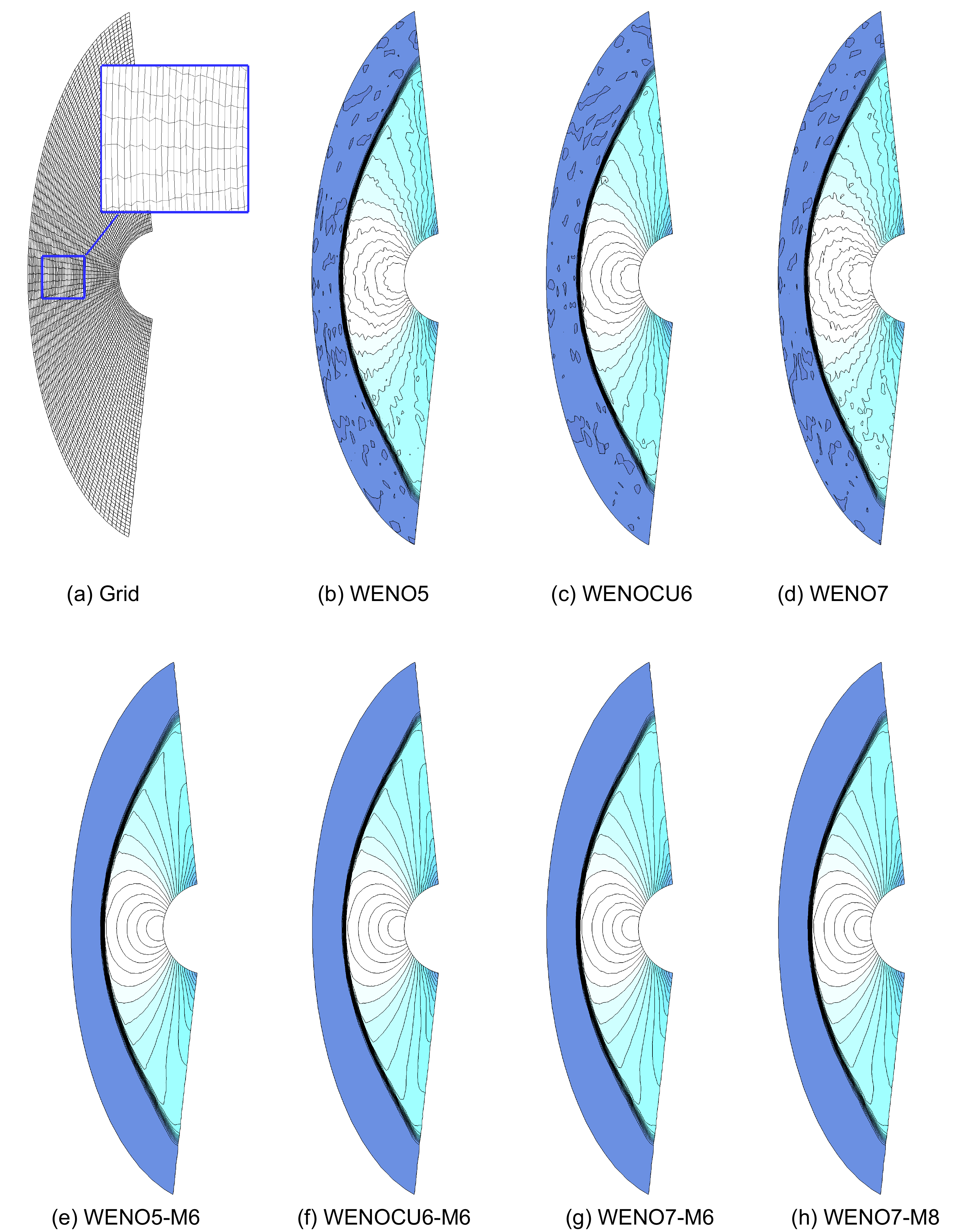}
    \caption{The pressure contours ranging from $1.2$ to $5.4$ of the supersonic flow past a cylinder.}
    \label{cylinder}
\end{figure}

The supersonic flow past a cylinder~\cite{jiang_efficient_1996} is solved to demonstrate the potential application of the present schemes on supersonic flows with curved walls. The initial condition is a uniform flow state with $Ma=2$.
The slip wall boundary condition is imposed at the cylinder surface, and the supersonic inflow and outflow boundary condition are employed at the left and right boundaries, respectively, with the grid generated by
\begin{equation}
\begin{aligned}
 x&=\left( R_x-(R_x-1) {\eta}'  \right) cos\left( \theta(2\xi'-1)\right),\\
 y&=\left( R_y-(R_y-1) {\eta}'  \right) sin\left( \theta(2\xi'-1) \right),\\
  \xi'&=\dfrac{\xi-1}{i_{max}-1}, \xi=i+0.2\phi_i , \\
 \eta'&=\dfrac{\eta-1}{j_{max}-1}, \eta=j+0.2\sqrt{1-\phi_i^2},
\end{aligned}
\end{equation}  
where $\theta=5\pi/12$, $R_x=3$, $R_y=6$, and $\phi_i \sim \mathcal{U}(0,1)$. The number of grid points is $i_{max}=61$ and $j_{max}=81$. The free-stream pressure and density are $p=1$ and $\rho=\gamma$, respectively. 
Our simulations are conducted until a steady-state is reached, with a time step of  $\Delta t=0.005$.
As shown in Fig. \ref{cylinder}(b-d), the grid randomness in Fig. \ref{cylinder}(a) leads to significant errors in the original WENO schemes, even in the smooth flow regions. However, our methods capture the bow shock sharply, and generate smooth contours away from the shock, without generating spurious free-streaming violations, as demonstrated by the pressure contours in Fig.~\ref{cylinder} (e-h).

\subsubsection{Transonic flow past a NACA0012 airfoil}\label{sec:NACA0012}
\begin{figure} [htbp]
    \center 
    \includegraphics[width=\textwidth]{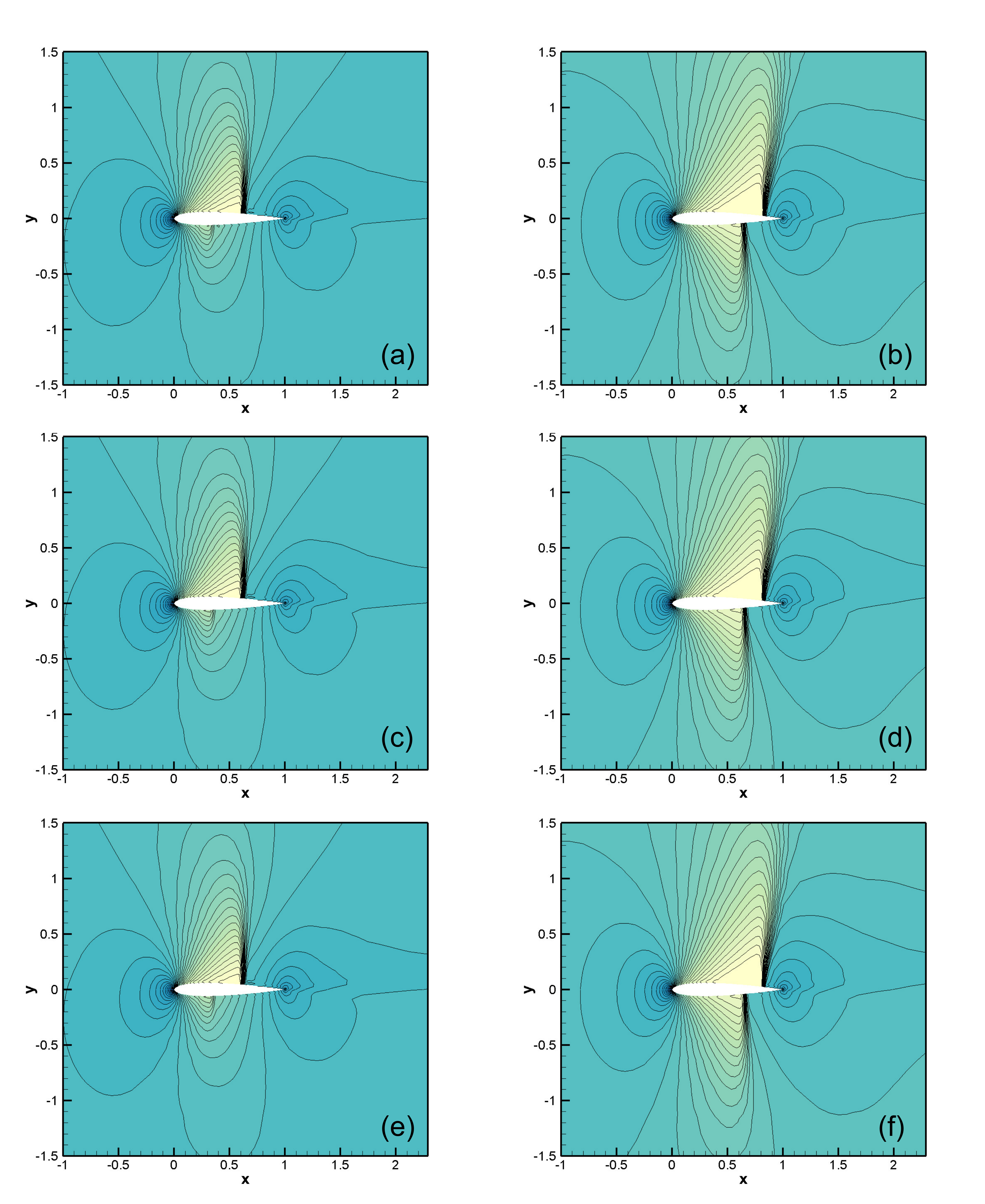}
    \caption{The Mach number contours ranging from $0.172$ to $1.325$ of NACA0012 airfoil for case 1 (the left column) and case 2 (the right column).}
    \label{naca0012_contours}
\end{figure}

\begin{figure} [htbp]
    \center 
    \includegraphics[width=\textwidth]{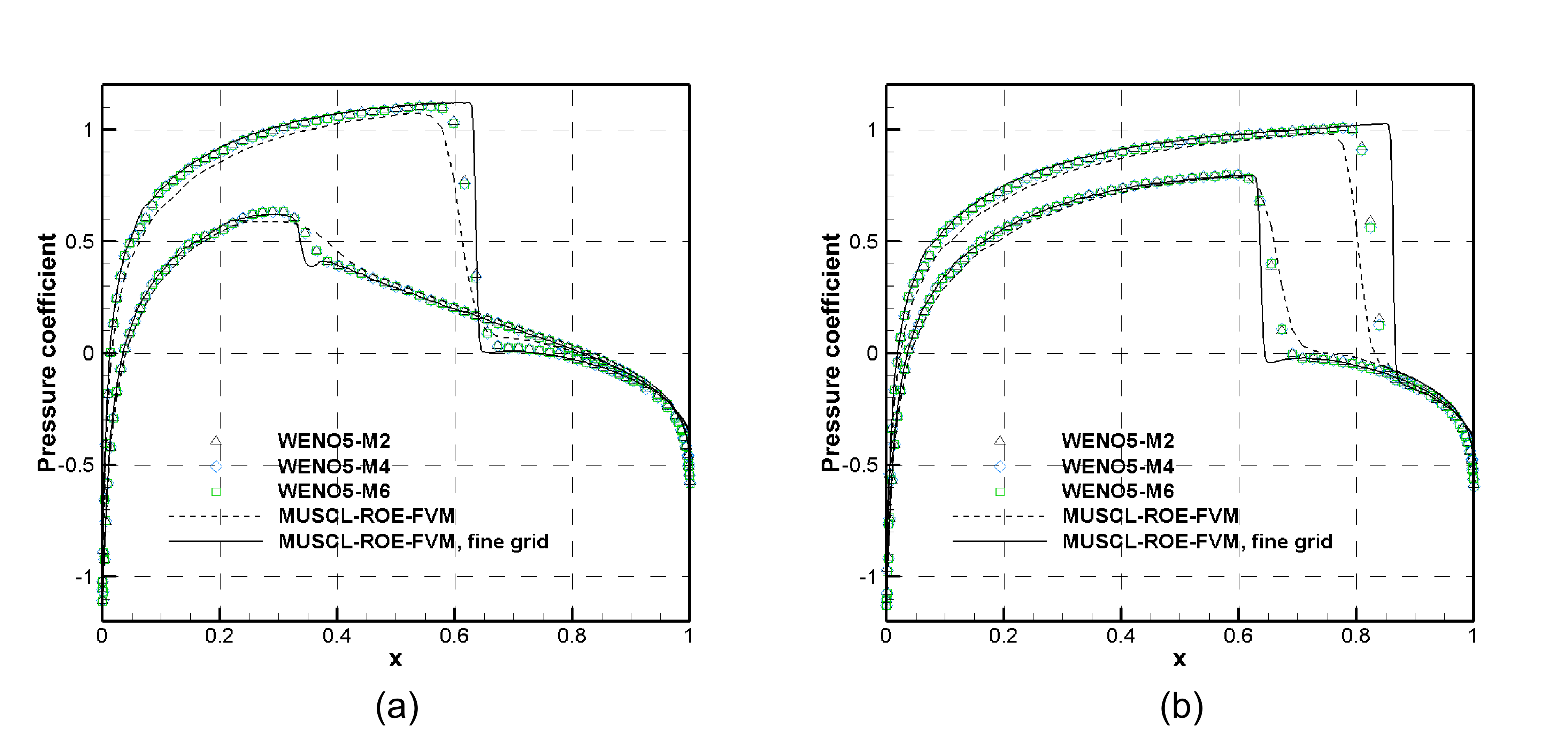}
    \caption{The pressure coefficient distributions on the surface of NACA0012 airfoil for case 1 (a) and case 2 (b). }
    \label{naca0012_Cp_cmp}
\end{figure}

Now we test the application of the proposed schemes on predicting the aerodynamic problems by considering the transonic flow past a NACA0012 airfoi. The simulations are performed on a coarse grids with a resolution of $160 \times 32$, with two setups: 1) Mach number $Ma=0.8$, angle of attack $AOA=1.25^{\circ}$ (case 1); 2) $Ma=0.85$, $AOA=1.0^{\circ}$ (case 2). The Mach-number contours in Fig. \ref{naca0012_contours} show that the WENO5-M2, WENO5-M4 and WENO-M6 schemes achieve similar flow structures in both two cases. And the corresponding pressure coefficient distributions along the airfoil surface of those schemes essential coincide with each other, see Fig. \ref{naca0012_Cp_cmp}. This observation indicate that the order of the metric discretization does not significantly affect the overall prediction accuracy on such a coarse resolution, which is in agreement with the conclusion in Sec. \ref{sec:Isotropic vortex}. Clearly, in comparison with the FVM based on 2nd-order MUSCL reconstuction and Roe scheme (hereafter referred to as MUSCL-ROE-FVM), improved agreement with the high-resolution ($1280 \times 177$ cells) reference MUSCL-ROE-FVM solution is achieved by the present WENO schmes, especially on the upper airfoil surface, with the shockwave being caputred more sharply.


\subsection{Vicous wall-bounded flows}
\subsubsection{Laminar vortex shedding from a circular cylinder} \label{sec:cylinder_unsteady}
\begin{figure} [htbp]
    \center 
    \includegraphics[width=0.9\textwidth]{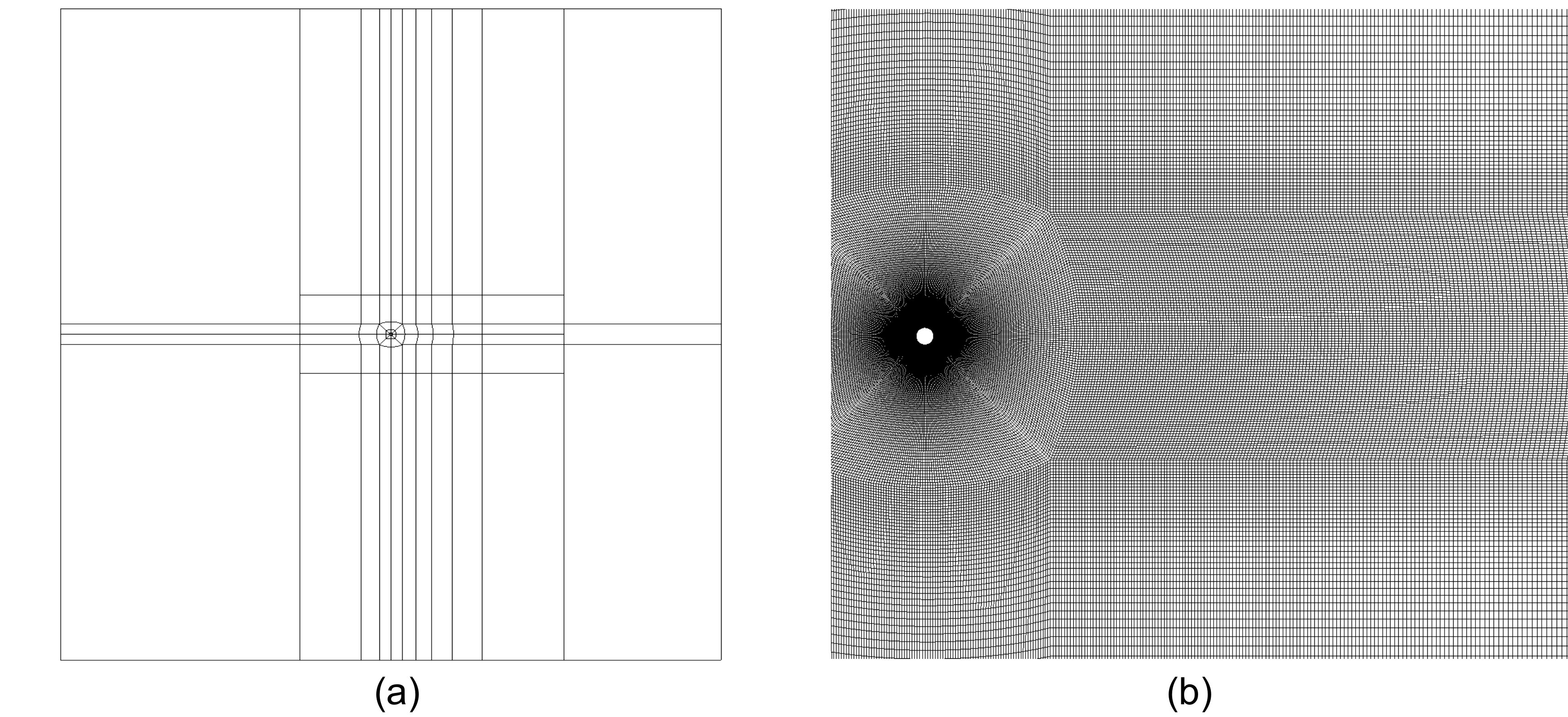}
    \caption{The block topology (a) and the enlarged view of the grid points (b) of the circular cylinder.}
    \label{cylinder_vortex_grid}
\end{figure}

\begin{figure} [htb]
    \center 
    \includegraphics[width=0.7\textwidth]{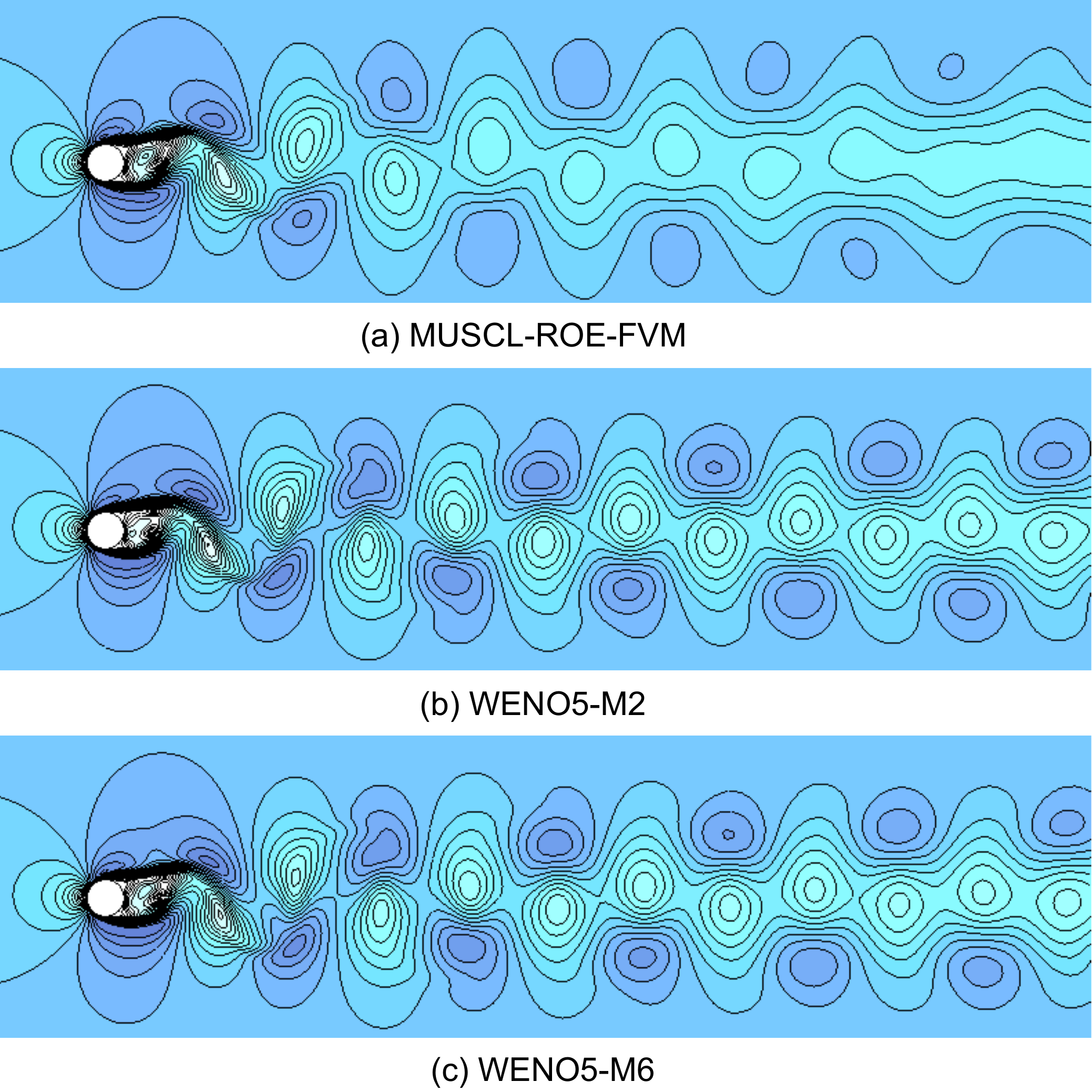}
    \caption{The Mach number contours ranging from $0$ to $0.26$ of the laminar vortex shedding from a circular cylinder at $Re_D=150$.}
    \label{cylinder_vortex}
\end{figure}

After testing the performce of the proposed method on inviscid flows, now we consider the 2D laminar subsonic flows past the circular cylinder at $Re_D=150$ \cite{Nichols2006Validation} to validate the accuracy of the FP schemes on resolving unsteady vortex shedding in viscous flows. The simulations are performed on a mesh with $177400$ cells, and a far-field boundary with $Ma=0.2$ and $T=277.78K$ is located $200$ diameters from the center of the cylinder, as illustrated in Fig. \ref{cylinder_vortex_grid}. To reduce the computational costs, the implicit dual-time-stepping algorithm with Lower-Upper Symmetric-Gauss-Seidel (LUSGS) \cite{1991Time, SEOKKWAN1988LUSGS} is employed to the simulations, with $\Delta t=0.001s$ for physical time marching and a CFL number of $20$ for sub-iteration
. In Fig. \ref{cylinder_vortex}, the Mach number contours indicate that the vortex shed from the cylinder in the wake is well resolved by the WENO5-M2 and WENO5-M6 schemes, and is significantly dissipated in the 2nd-order MUSCL-ROE-FVM result.
As listed in Table~\ref{cylinder_vortex_table}, compared to experimental data in Ref. \cite{Nichols2006Validation}, both the predicted averaged drag coefficient $C_D$ and lift Strouhal number $St$ of our schemes show smaller error than of the MUSCL-ROE-FVM. In addition, the WENO5-M2 scheme achieves approximately the same results with those in the WENO5-M6 scheme, both the vortical stuctures in the wake and statistic data ($C_D$ and $St$), which indicates the order of metrics and Jacobian has no significant influences on the under-resolved simulations, as in Sec. \ref{sec:Isotropic vortex} and Sec. \ref{sec:NACA0012}.


\begin{table}[htbp]	
   \centering
	\caption{The experimental and simulated averaged drag coefficient ($C_D$) and lift Strouhal number ($St$) of the laminar vortex shedding from a circular cylinder.}	 
	
	\begin{tabular}{ccccc}
  
   \hline    
         &Exp.     &MUSCL-ROE-FVM &WENO5-M2 &WENO5-M6 \\   
         \hline  
   $C_D$ &1.34     &1.312          &1.334    &1.334  \\
   \hline  
 $St$ &0.179-0.182 &0.172           &0.182     &0.182   \\
		\hline  
	\end{tabular}
	\label{cylinder_vortex_table}
\end{table}

\subsubsection{Subsonic flow past a three-element airfoil (30P30N)}\label{sec:sphere}
\begin{figure} [htbp]
    \center 
    \includegraphics[width=\textwidth]{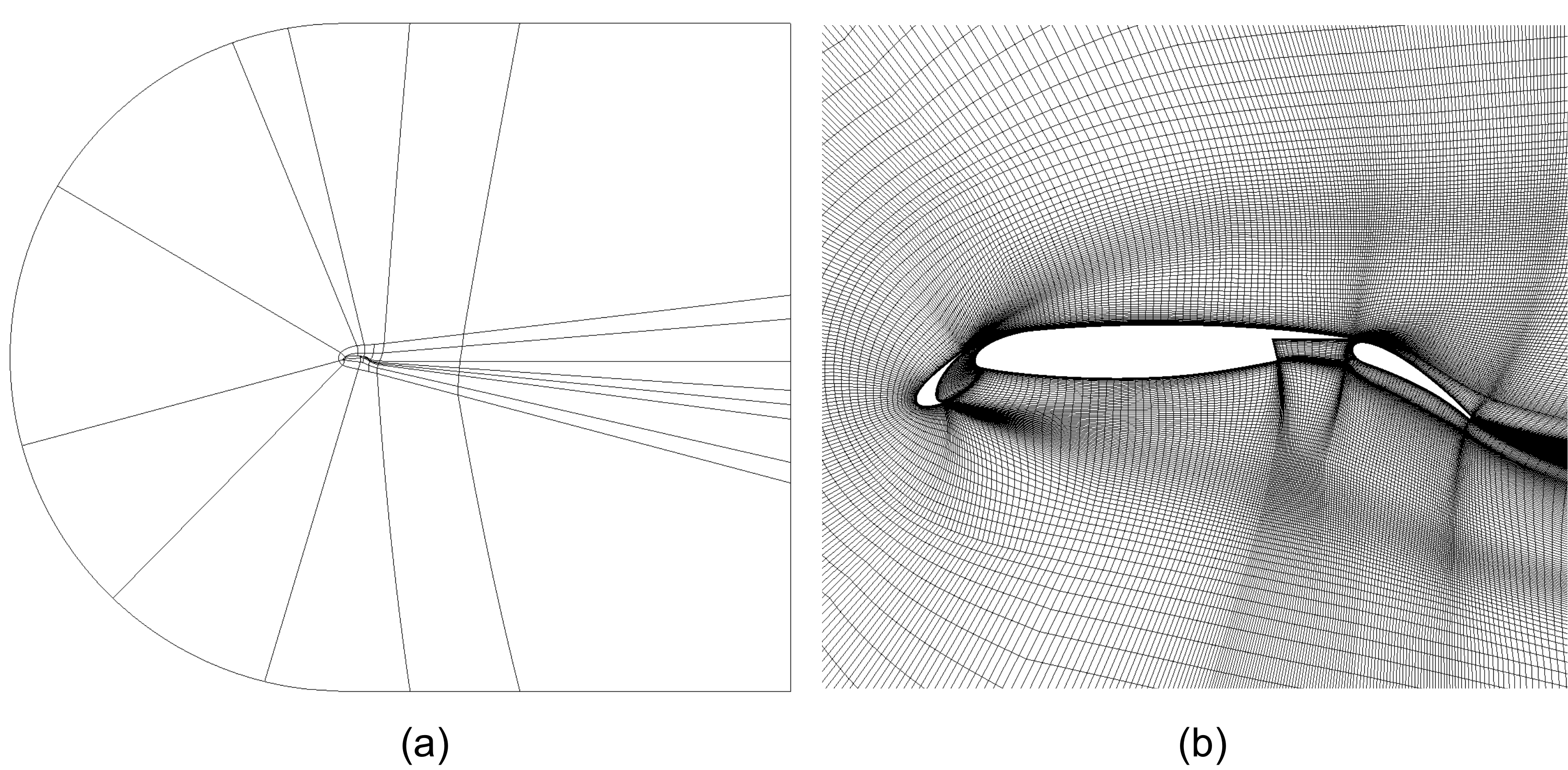}
    \caption{The block topology (a) and the enlarged view of the grid points (b) of the 30P30N three-element airfoil.}
    \label{30P30N_grid}
\end{figure}

In this section, the subsonic flow past a three-element airfoil (30P30N) is simulated to demonstrate the capability and accuracy of the present FP schemes on simualting  aerodynamic problems of complex geometries on turbulent conditions.
Here, the far-field boundary condition with $Ma=0.2$ and $AOA=8^{\circ}$ is imposed on the boundary of the computational domain which is discretized by a mesh containing $10^5$ cells, as shown in Fig. \ref{30P30N_grid}. The steay solution is obtained by the implicit time marching with LUSGS and the SST turbulence model \cite{Menter1994Two}.

As illustrated in Fig. \ref{30P30N_grid}(a-d), WENO5-M2, WENO5-M6
and Linear-Upwind5-M6 produce similar continuous and smooth Mach-number contours around the 30P30N airfoil. The slip flows at the tail of the slat is thinner in these FP schemes, and at the flap tail we observe vortex shedding which is absent in the 2nd-order MUSCL-ROE-FVM. In addition, as in Fig. \ref{30P30N_grid}(e), the pressure coefficients distributions along the airfoil surface in the present schemes match the experimental data in Ref. \cite{Cai2006A} very well, while are significantly underestimated on the upper surface in the 2nd-order MUSCL-ROE-FVM results.

\begin{figure} [p]
    \center 
    \includegraphics[width=0.8\textwidth]{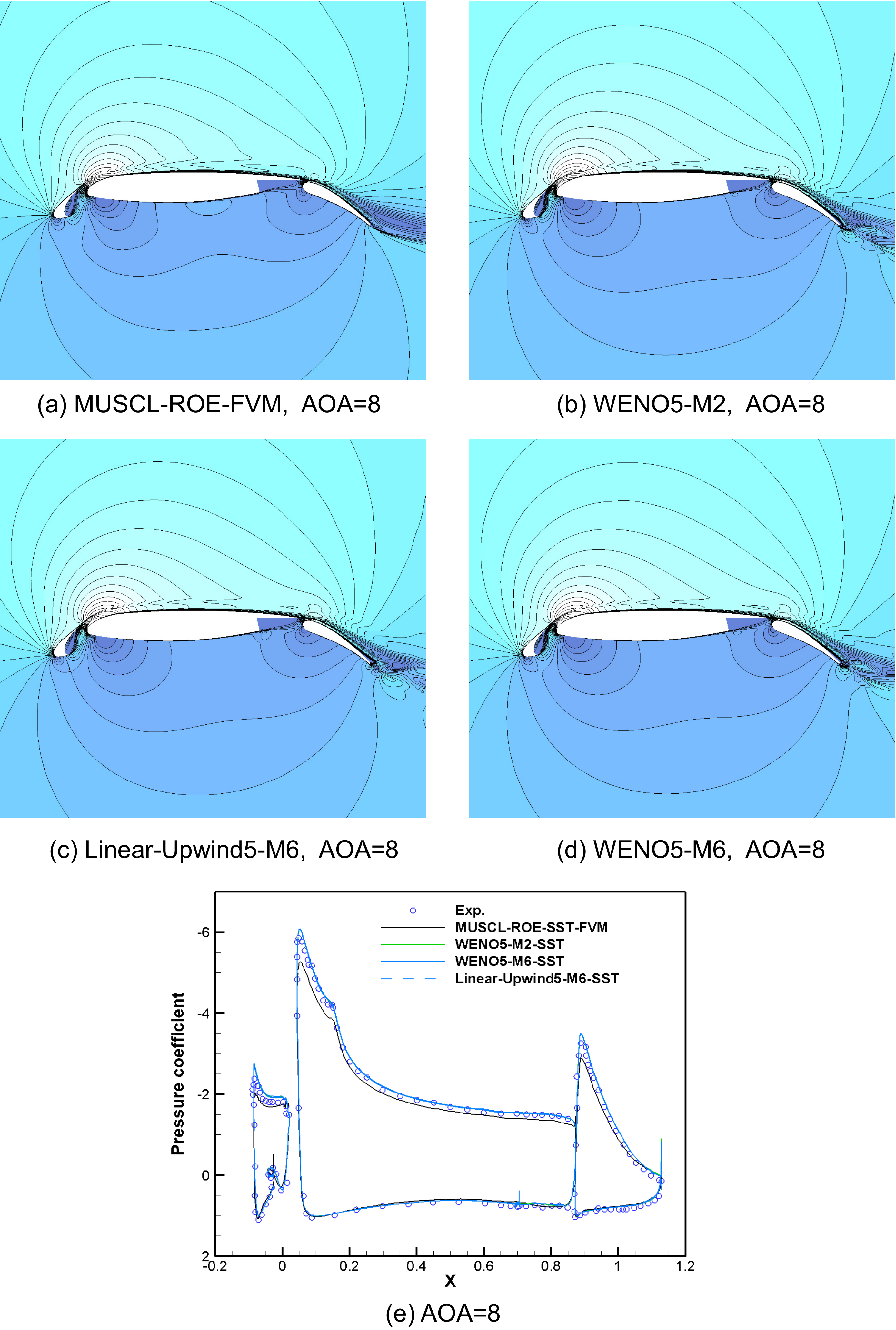}
    \caption{The Mach number contours ranging from $0.05$ to $0.5$ in different schemes (a-d) and the pressure coefficient distributions on the surface of 30P30N airfoil (e).}
    \label{30P30N_contours}
\end{figure}

\subsubsection{Transonic flow pass the ONERA M6 wing}\label{sec:M6_wing}
In addtion to the 2D test cases above, we also select the tansonic flow around the ONERA M6 wing to validate the robustness of our schemes on predicting turbulent flows in three dimensions. We perform the simulations under a condition of $Ma=0.84$, and $AOA=3.06^{\circ}$, with LUSGS for time advancing and a mesh of $294912$ cells for space discretization. The Reynolds number of $Re_l=1.172 \times 10^7$ is sufficiently large that the SST turbulence model is invoked, where $l=0.64607m$ is the mean aerodynamic chord.

\begin{figure} [p]
    \center 
    \includegraphics[width=\textwidth]{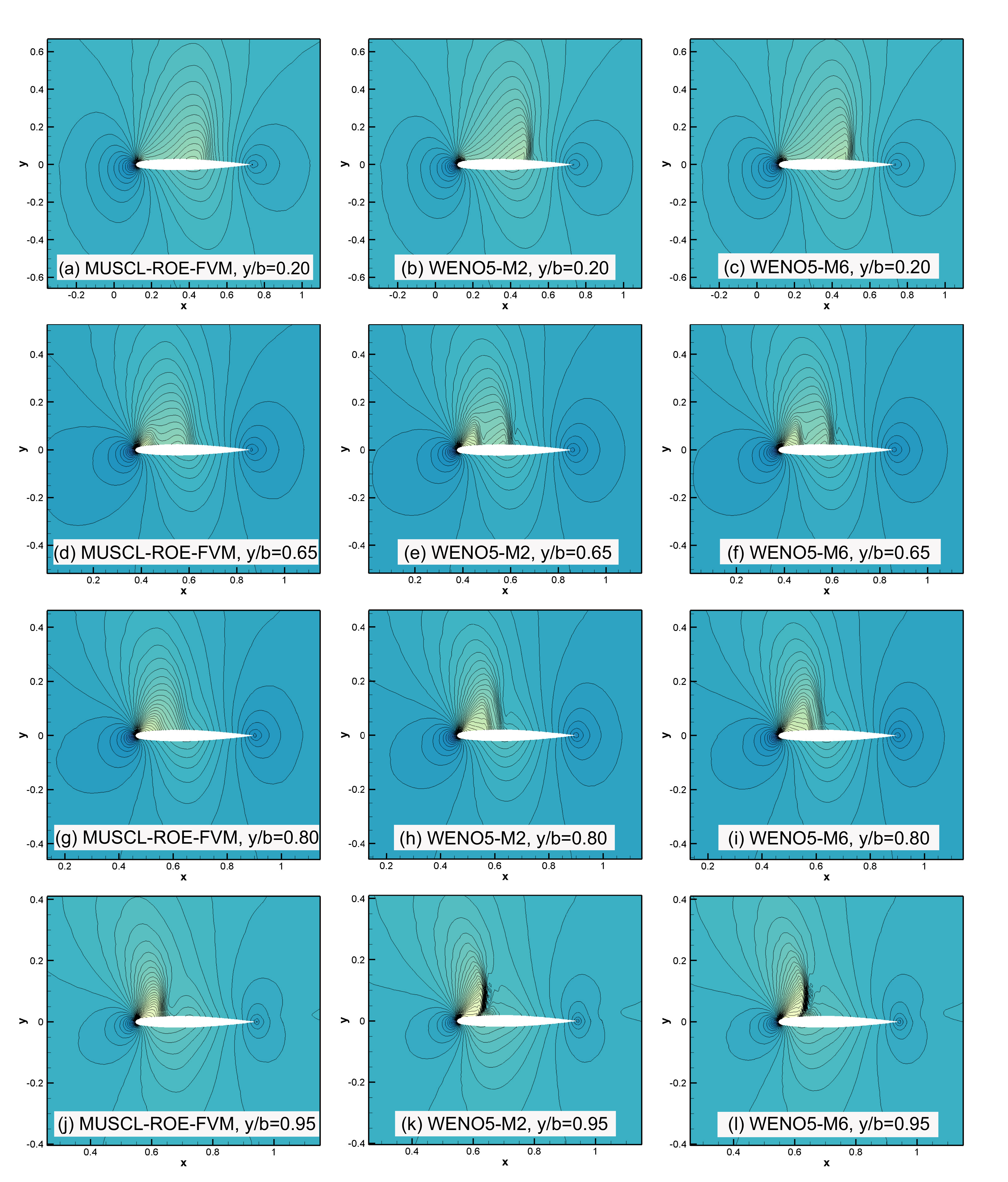}
    \caption{The pressure contours ranging from $1.2 \times 10^5$Pa to $4.6 \times 10^5$Pa of the ONERA M6 wing at four different spanwise cross-sections.}
    \label{M6_contours}
\end{figure}

\begin{figure} [p]
    \center 
    \includegraphics[width=\textwidth]{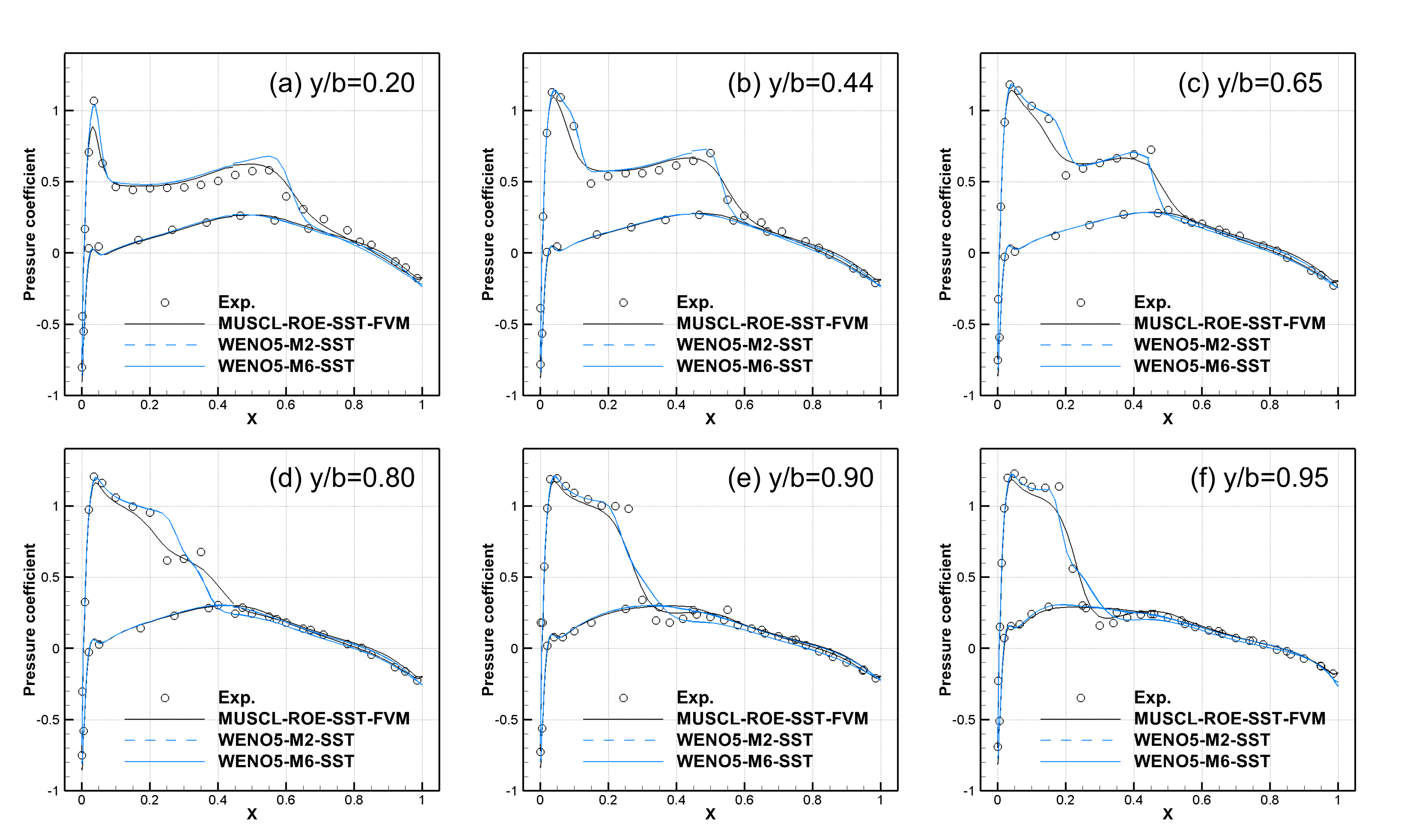}
    \caption{The pressure coefficient distributions at $6$ cross-sections along the spanwise direction of the ONERA M6 wing surface.}
    \label{M6_Cp}
\end{figure}

The pressure contours in Fig. \ref{M6_contours} indicate that both the WENO5-M2 and WENO5-M6 schemes obtain a better result than MUSCL-ROE-FVM, as the flow filed near the shockwave is more clearly captured. As illustrated in Fig. \ref{M6_Cp}, compared to MUSCL-ROE-FVM, the pressure distributions of WENO5-M2 and WENO5-M6 are in a better agreement with the experimental results \cite{Schmitt1979PressureDO} at the $6$ selected cross-sections along the spanwise direction. Again, we emphasize that the order of metric does not actually affect the predicting accuracy in this 3D problem, as demonstraded by the qualitative comparision in Fig. \ref{M6_contours} and the quantitative comparision in Fig. \ref{M6_Cp}, between WENO5-M2 and WENO5-M6.

\section{Conclusion} \label{sec:conlusion}
The proposed free-stream preserving method in this work suggest to reformulate the dissipation
of the upwind conservative finite-difference schemes, in order to eliminate the free-stream preserving violation errors on curvilinear grids. In particular, to construct the reformulated dissipation, the reference variables are subtracted from the cell-center conservative variables and fluxes for the linear upwind schemes and the WENO schemes, respectively. In such a way, the proposed numerical approach and its algorithmic formulation have the following main properties: (1) no modification in this method is required on the original formulation of the upwind conservative finite-difference schemes; (2) this method is more flexible in practical engineering problems, as the metrics and fluxes can be discretized in an inconsistent manner, with the free-steaming condition preserved.  
The simple test cases demonstrate that the proposed method strictly preserve the geometric conservation law on stationary curvilinear grids, with a formal high-order accuracy. The accuracy and robustness of this method are  assessed in various fluid dynamics problems involving complex geometrics, e.g. the NACA0012 airfoil, the 30P30N airfoil, and the ONERA M6 wing, where the shockwaves, vortical structures, and turbulent flows are well predicted. In addition, it can be concluded from a range of complex validation cases that the accuracy of metrics discretization show no significant effect on under-resolved numerical simulations, which are usually encountered in piratical problems where low-quality and coarse-resolution grids are widely used. 

\section*{Appendix}\label{appd}
The final cell-face fluxes reconstructed by the proposed WENO scheme can be expressed by the summation of the central and dissipation part~\cite{nonomura_new_2015, zhu_free-stream_2019} ,
\begin{equation}\label{eq:WENO_cen_dis}
\begin{aligned}
\boldsymbol{ \tilde{F} }_{i+1/2}^{'}=&\boldsymbol{ \tilde{F} }_{i+1/2}^{'+}+\boldsymbol{ \tilde{F} }_{i+1/2}^{'-}\\
=&\sum\limits_s \boldsymbol{R}_{i+1/2}^s {f}_{i+1/2}^{s,+}  +  \sum\limits_s \boldsymbol{R}_{i+1/2}^s {f}_{i+1/2}^{s,-}\\
=&\dfrac{1}{60} \left( \boldsymbol{\tilde{F}}_{i-2} -8 \boldsymbol{\tilde{F}}_{i-1}+37 \boldsymbol{\tilde{F}}_{i} +37\boldsymbol{\tilde{F}}_{i+1}-8 \boldsymbol{\tilde{F}}_{i+2}+\boldsymbol{\tilde{F}}_{i+3} \right)-\boldsymbol{\hat{F}}_{i+1/2}\\
&  -\dfrac{1}{60} \sum\limits_s \boldsymbol{R}_{i+1/2}^s \left[ \left(20\omega_0^+-1 \right) \hat{f}_{i,1}^{s,+}-\left( 10\omega_0^+ +10\omega_1^+ -5 \right)\hat{f}_{i,2}^{s,+} +\hat{f}_{i,3}^{s,+}\right]\\
&  +\dfrac{1}{60} \sum\limits_s \boldsymbol{R}_{i+1/2}^s \left[ \left(20\omega_0^- -1 \right) \hat{f}_{i,1}^{s,-}-\left( 10\omega_0^- +10\omega_1^- -5 \right)\hat{f}_{i,2}^{s,-} +\hat{f}_{i,3}^{s,-}\right],
\end{aligned}
\end{equation}
where
\begin{equation}\label{eq:WENO_hatF}
\begin{aligned}
\boldsymbol{\hat{F}}_{i+1/2}&=\dfrac{1}{60} \left[ \left( \boldsymbol{\tilde{F}}_{i-2}^{*}-\boldsymbol{\tilde{F}}_{i+1/2}^{*} \right) -8 \left( \boldsymbol{\tilde{F}}_{i-1}^{*}-\boldsymbol{\tilde{F}}_{i+1/2}^{*} \right)+37 \left( \boldsymbol{\tilde{F}}_{i}^{*}-\boldsymbol{\tilde{F}}_{i+1/2}^{*} \right) \right. \\
& \quad \left. +37\left( \boldsymbol{\tilde{F}}_{i+1}^{*}-\boldsymbol{\tilde{F}}_{i+1/2}^{*} \right)-8 \left( \boldsymbol{\tilde{F}}_{i+2}^{*}-\boldsymbol{\tilde{F}}_{i+1/2}^{*} \right)+\left( \boldsymbol{\tilde{F}}_{i+3}^{*}-\boldsymbol{\tilde{F}}_{i+1/2}^{*} \right) \right]\\
&=\dfrac{1}{60}  \left( \boldsymbol{\tilde{F}}_{i-2}^{*} -8 \boldsymbol{\tilde{F}}_{i-1}^{*}+37 \boldsymbol{\tilde{F}}_{i}^{*} +37 \boldsymbol{\tilde{F}}_{i+1}^{*}-8  \boldsymbol{\tilde{F}}_{i+2}^{*}+ \boldsymbol{\tilde{F}}_{i+3}^{*}\right) -\boldsymbol{\tilde{F}}_{i+1/2}^{*}
\end{aligned}
\end{equation}
and
\begin{equation} \label{eq:WENO_cen_dis_sub1}
\begin{aligned}
\hat{f}_{i,r+1}^{s,+} =&\widetilde{f}_{i+r+1}^{s,+}-3\widetilde{f}_{i+r}^{s,+}+3\widetilde{f}_{i+r-1}^{s,+}-\widetilde{f}_{i+r-2}^{s,+} , \qquad r=0,1,2\\
=&\dfrac{1}{2} \boldsymbol{L}_{i+1/2}^s \left(  \boldsymbol{\tilde{F}}_{i+r+1}-3\boldsymbol{\tilde{F}}_{i+r}+3\boldsymbol{\tilde{F}}_{i+r-1}-\boldsymbol{\tilde{F}}_{i+r-2} \right) \\
&+\dfrac{1}{2}\lambda^s \boldsymbol{L}_{i+1/2}^s\left( \boldsymbol{\tilde{Q}}_{i+r+1}-3\boldsymbol{\tilde{Q}}_{i+r}+3\boldsymbol{\tilde{Q}}_{i+r-1}-\boldsymbol{\tilde{Q}}_{i+r-2} \right)\\
&-\dfrac{1}{2} \boldsymbol{L}_{i+1/2}^s \left(  \boldsymbol{\tilde{F}}_{i+r+1}^{*}-3\boldsymbol{\tilde{F}}_{i+r}^{*}+3\boldsymbol{\tilde{F}}_{i+r-1}^{*}-\boldsymbol{\tilde{F}}_{i+r-2}^{*} \right)\\
&-\dfrac{1}{2}\lambda^s \boldsymbol{L}_{i+1/2}^s\left( \boldsymbol{\tilde{Q}}_{i+r+1}^{*}-3\boldsymbol{\tilde{Q}}_{i+r}^{*}+3\boldsymbol{\tilde{Q}}_{i+r-1}^{*}-\boldsymbol{\tilde{Q}}_{i+r-2}^{*} \right),
\end{aligned}
\end{equation}

Finally, another form of the free-stream preserving linear upwind scheme can be obtained by imposing the optimal weights to the WENO scheme, as shown in Eq.~\eqref{eq:new_form_linear_upwind}.

Using the optimal weights in Eq.~\eqref{eq:weno_wights}, Eq. \eqref{eq:WENO_cen_dis} recovers to the linear upwind FP scheme,
\begin{equation}\label{eq:new_form_linear_upwind}
\begin{aligned}
\boldsymbol{ \tilde{F} }_{i+1/2}^{'}
=&\dfrac{1}{60} \left( \boldsymbol{\tilde{F}}_{i-2} -8 \boldsymbol{\tilde{F}}_{i-1}+37 \boldsymbol{\tilde{F}}_{i} +37\boldsymbol{\tilde{F}}_{i+1}-8 \boldsymbol{\tilde{F}}_{i+2}+\boldsymbol{\tilde{F}}_{i+3} \right)-\boldsymbol{\hat{F}}_{i+1/2}\\
&+\dfrac{1}{60}\sum\limits_s \boldsymbol{R}_{i+1/2}^s  \lambda^s \boldsymbol{L}_{i+1/2}^s \left[ \left(\boldsymbol{\tilde{Q}}_{i-2}-\boldsymbol{\tilde{Q}}_{i-2}^*\right)-5\left(\boldsymbol{\tilde{Q}}_{i-1}-\boldsymbol{\tilde{Q}}_{i-1}^*\right)+10\left(\boldsymbol{\tilde{Q}}_{i}-\boldsymbol{\tilde{Q}}_{i}^*\right) \right.\\ 
 &  \left. -10\left(\boldsymbol{\tilde{Q}}_{i+1}-\boldsymbol{\tilde{Q}}_{i+1}^*\right)+5\left(\boldsymbol{\tilde{Q}}_{i+2}-\boldsymbol{\tilde{Q}}_{i+2}^*\right)-\left(\boldsymbol{\tilde{Q}}_{i+3}-\boldsymbol{\tilde{Q}}_{i+3}^*\right) \right],
\end{aligned}
\end{equation}
which clearly is able to use inconsistent discretization methods for the metrics and fluxes, unlike the linear upwind FP scheme in Eq. \eqref{eq:final_form_linear_upwind} where the sufficient conditon of Deng and Abe ~\cite{ DENG201390,abe2014geometric} is required.
\section*{Acknowledgements}
This work was supported by the National Natural Science Foundation of China (Grant Nos. 11902271 and 91952203), the Fundamental Research Funds for the Central Universities of China (No. G2019KY05102), the Foundation of National Key Laboratory (No. 6142201190303), and
111 project on ``Aircraft Complex Flows and the Control'' (Grant No. B17037). The authors gratefully acknowledge Yujie Zhu for discussion on the simulation setups.

\section*{References}

\bibliographystyle{elsarticle-num}
\bibliography{ref}

\end{document}